\documentclass[12pt]{article}
\usepackage{amsmath,amsthm,amsfonts,amssymb}

\title{The Tits Alternative for $Out(F_n)$ II:\\ A Kolchin Type Theorem}
\author{Mladen Bestvina, Mark Feighn, and Michael Handel}
\date{\today\\ preliminary version}

\newenvironment{pf}{\begin{proof}}{\end{proof}}
\def\ffs{free factor system}
\def\xx{{\cal X}}
\def\oa{{\cal O}}
\def\ttt{G}
\def\tmm{f}

\newtheorem{thm}{Theorem}[section]
\newtheorem{lem}[thm]{Lemma}
\newtheorem{cor}[thm]{Corollary}
\newtheorem{prop}[thm]{Proposition}

\newtheorem*{bouncing kolchin}{Theorem \ref{bouncing ends well}}{}
\newtheorem*{no number thm}{Theorem}{}
\newtheorem*{weak tits}{Theorem (Tits Alternative for
$Out(\f)$)}{}
\newtheorem*{graph kolchin macro}{Theorem \ref{graph kolchin} (Kolchin Theorem for $Out(F_n)$}{}
\newtheorem*{tree kolchin}{Theorem \ref{kolchin for trees}}{}

\newtheorem{sublem}[thm]{Sublemma}
\newtheorem{prop-def}[thm]{Definition-Proposition}

\newtheorem{stoa}{Theorem (Solvable implies abelian)}

\theoremstyle{definition}
\newtheorem{defn}[thm]{Definition}

\newtheorem{notn}[thm]{Notation}
\newtheorem{com}[thm]{Remark}
\newtheorem{ex}[thm]{Example}

\newtheorem{tits}{Theorem (The Tits Alternative for $Out(F_n)$)}

\theoremstyle{remark}

\begin{document}


\maketitle


\tableofcontents

\def\Cal{\cal}

\newcommand{\h}{\frak H} 
\newcommand{\x}{{\Cal X}_S} 
\newcommand{\ve}{{\Cal V}} 
\newcommand{\PG}{PG_0(1)}

\newcommand{\vs}{\x} 
\newcommand{\tbar}{\overline T} 
\newcommand{\e}{\epsilon}

\def\G{G}
\def\GG{{\cal H}}
\def\H{{\cal H}}
\def\R{\cal R}
\def\f{{F_n}} 
\def\U{UPG_{\f}}
\def\Beta{\cal B}
\def\pgif{$PGIF$}

\newcommand{\n}{n} 
\newcommand{\bd}{\partial} \bibliographystyle{amsalpha}


\def\hom(#1){H_1(#1,\mathbb Z/3\mathbb Z)}
\def\l(#1,#2){\ell_{#1}(#2)}
\def\roster{\begin{enumerate}}
\def\endroster{\end{enumerate}}
\def\definition{\begin{defn}}
\def\enddefinition{\end{defn}}
\def\subhead{\subsection\{}
\def\endsubhead{\}}
\def\head{\section\{}
\def\endhead{\}}
\def\example{\begin{ex}}
\def\endexample{\end{ex}}
\def\ves{\vs}
\def\cv{CV}

\def\Q{{\cal Q}}
\def\q{f}
\def\r{g}
\def\Z{{\cal Z}}
\def\S{{\cal S}}

\section{Introduction and Outline} \label{intro}

 Recent
years have seen a development of the  theory for
$Out(\f)$, the outer automorphism group of the free group $\f$
of rank
$n$, that is modelled on Nielsen-Thurston theory for surface
homeomorphisms.
As mapping classes have either  exponential or linear growth
rates, so  free
group outer automorphisms have either  exponential or
polynomial growth rates. (The degree of the polynomial can be
any integer
between 1 and $n-1$, see
\cite{bh:tracks}.) In \cite{bfh:tits1} we considered individual
automorphisms, with primary emphasis on those with exponential
growth rates.  In this paper we focus on subgroups  of $Out(\f)$,
all of
whose elements have  polynomial growth rates.

     To remove certain technicalities arising from finite order
phenomena, we
restrict our attention to  those    polynomially growing outer
automorphisms $\oa$ whose induced automorphism of
$H_1(\f;\mathbb Z)\cong {\mathbb Z}^n$ is unipotent.  We say that
such an outer
automorphism is  {\it unipotent}; we also say that $\oa$ is
a $UPG(\f)$ (or just a $UPG$) outer automorphism. A subgroup
of $Out(\f)$ is
called a {\it $UPG$ subgroup} if each element is  $UPG$.
We prove (Proposition \ref{unichar}) that any polynomially
growing outer
automorphism that acts trivially in ${\mathbb Z}/3\mathbb Z$-homology is
unipotent.  Thus every subgroup of polynomially growing outer
automorphisms has
a finite index $UPG$ subgroup.

      The archetype for the main theorem of this paper comes from
linear
groups.  A linear map is unipotent if and only if it  has a basis
with respect to
which it is upper triangular with 1's on the diagonal. A celebrated
theorem of
Kolchin \cite{se:lie} states that for any group of unipotent
linear maps there is
a basis with respect to which  all elements of the group are upper
triangular with 1's on the diagonal.

    There is an analagous result for mapping class groups.  We say
that a mapping class is unipotent if it has linear growth and if the
induced linear map on first homology is unipotent. The Thurston
classification theorem implies that a mapping class is unipotent if
and only if it is represented by a composition of Dehn twists in
disjoint simple closed curves.  Moreover, if a pair of unipotent
mapping classes belong to a unipotent subgroup, then their twisting
curves can not have transverse intersections (see for example
\cite{blm:s2a}).  Thus every unipotent mapping class subgroup has a
characteristic set of disjoint simple closed curves and each element
of the subgroup is a composition of Dehn twists along these
curves.  As in the linear case, in which the basis does not depend on
the individual linear maps in a the unipotent subgroup, here the
twisting curves do not depend on the individual mapping classes.

     Our main theorem is the analogue of Kolchin's theorem
for $Out(\f)$. Recall \cite{cv:moduli} that a marked graph is a
graph (1-dimensional
CW-complex) equipped with a homotopy equivalence from the
rose with
$n$ petals (whose fundamental group is permanently identified
with
$F_n$). A homotopy equivance $\tmm:\ttt\to\ttt$ on a marked
graph $G$ induces
an outer automorphism of the fundamental group of $G$ and
therefore an
element $\oa$ of $Out(\f)$; we say that $\tmm:\ttt\to\ttt$ is a
{\it representative} of $\oa$.

     Suppose that $G$ is a marked graph and
that $\emptyset=\ttt_0\subset\ttt_1\subset\cdots\subset
\ttt_K=\ttt$  is a filtration of $G$ where $\ttt_{i}$ is
obtained  from
$\ttt_{i-1}$ by adding a single edge $E_i$.  A homotopy
equivalence
$\tmm:\ttt\to\ttt$ is  {\it upper
triangular with respect to the filtration} if each
$\tmm(E_i)=v_iE_iu_i$
where $u_i$ and $v_i $ are loops in $\ttt_{i-1}$. If the choice of
filtration is
clear then we  simply say that $\tmm:\ttt\to\ttt$ is   {\it upper
triangular}. We refer to the $u_i$'s and $v_i$'s as  {\it
suffixes} and  {\it prefixes} respectively.

     An outer automorphism is $UPG$ if and only if it has a
representative
that is upper triangular with respect to some filtered marked graph $G$
(see Section \ref{unipotent}).

     For any filtered marked graph $G$ let $\Q$ be
the set of  upper triangular  homotopy equivalences of $G$ up to
homotopy
relative to the vertices of $G$.
By Lemma \ref{gp str}, $\Q$ is a group under the operation induced by composition.
    There is a  natural map from $\Q$ to $UPG(\f)$.  We say that
a subgroup of
$UPG(\f)$ is {\it filtered} if it lifts to a subgroup of $\Q$ for
some
filtered marked graph.
   We can now state our main theorem.

\begin{thm} {\rm\bf (Kolchin Theorem for $Out(F_n)$).} 
\label{graph kolchin} Every
finitely generated $UPG$ subgroup $\H$ of $Out(\f)$ is filtered.  The
number of edges of the filtered marked graph can be taken to be
bounded by $\frac{3n}{2}-1$ for $n>1$.\end{thm}

     It is an interesting question whether or not
the requirement that $\H$ be finitely generated is necessary or
just
an artifact of our proof.
\vskip .2cm
\noindent {\bf Question:} Is every $UPG$ group in $Out(\f)$
contained
in a finitely generated $UPG$ group?
\vskip .2cm

\begin{com}  \label{complicated} In contrast to unipotent
mapping class
subgroups, which are all finitely generated and abelian, $UPG$
subgroups of
$Out(\f)$ can be quite large. For example,  if $G$ is the rose on
$n$ petals, then a filtration on $G$ corresponds to an ordered
basis
$x_1,\cdots,x_n$  of $F_n$ and elements of $\Q$ correspond to
automorphisms of the form $x_i \mapsto a_ix_ib_i$ with
$a_i,b_i\in \langle
x_1,\cdots,x_{i-1}\rangle$. When $n > 2$, the image of $\Q$ in
$UPG(\f)$
contains $F_2\times F_2$.
\end{com}

This is the second  of two papers in which we establish
the Tits Alternative for $Out(\f)$.

\begin{tits} \label{tits theorem}
Let $\GG$ be any subgroup of $Out(\f)$. Then
either $\GG$ is virtually solvable, or
contains $F_2$.
\end{tits}

 For a proof of a special (generic) case, see \cite{bfh:tits0}.  The
relation between Theorem \ref{tits theorem} and Theorem \ref{graph
kolchin} is captured by the following corollary.

\begin{cor} \label{upg tits}
Every $UPG$ group $\H$ either contains $F_2$ or  is
solvable.
\end{cor}

\begin{pf} First assume that $\H$ is finitely generated. By
Theorem~\ref{graph kolchin} there is a marked graph $G$, a filtration
$\cal F$ and a subgroup $\Z$ of $\Q$ that projects isomorphically onto
$\H$. Let $i \ge 0$ be the largest parameter value for which every
element of $\Z$ restricts to the identity on $\ttt_{i-1}$. If $i =
K+1$, then $\Z$ is the trivial group and we are done.  Supppose then
that $i
\le K$.
By construction, each element of $\Z$ satisfies  $E_i
\mapsto v_iE_iu_i$ where
$v_i$ and $u_i$ are paths (that depend on the element of $\Z$) in
$\ttt_{i-1}$ and are therefore fixed by every element of $\Z$.  The
suffix map $\S : \Z \to \f$, which assigns the suffix $u_i$ to the
element of $\Z$, is therefore a homomorphism. If the image of $S$
contains $F_2$, then $\Z$ contains $F_2$ and we are done. If the image
of $S$ has rank one, then it can be identified with $\mathbb Z$ and
there is no loss in replacing $\Z$ with the kernel of $S$. If the
image of $S$ has rank zero, then $\Z$ is the kernel of $\S$. A similar
argument using prefixes instead of suffixes, allows us to replace $\Z$
with the subgroup of $\Z$ that has no non-trivial prefixes or suffixes
for $E_i$ and so restricts to the identity on $\ttt_i$. Upward
induction on $i$ now completes the proof when $\H$ is finitely generated. 
In fact, this argument shows that $\H$ is polycyclic and that the length
of the derived series is bounded by $\frac{3n}{2}-1$ for $n>1$. 

When $\H$ is not finitely generated, it can be represented as the
increasing union of finitely generated subgroups. If one of these
subgroups contains $F_2$, then so does $\H$, and if not then $\H$ is
solvable with the length of the derived series bounded by
$\frac{3n}{2}-1$.\end{pf}

\vskip.2cm

\noindent{\bf Proof of the Tits Alternative for $Out(F_n)$.} Theorem 1.3 of
 \cite{bfh:tits1} asserts that if
$\GG$ does not contain $F_2$, then there is an exact sequence
$$1\to
\H_0\to
\GG\to
{\cal A}\to 1$$ with $\cal A$ a finitely generated free abelian group
and with all elements of $\H_0$ of polynomial growth.

By passing to a subgroup of $\GG$
of finite index that acts trivially in ${\mathbb
Z}/3\mathbb Z$-homology, we may assume that $\H_0$ is a $UPG$ group
(see Proposition
\ref{unichar}). Since $\H_0$ does not contain
$F_2$, by Corollary \ref{upg tits}, $\H_0$ is solvable, and
thus
$\GG$ is also solvable.
\qed

\vskip .2cm

In \cite{bfh:tits3} we strengthen
the Tits
Alternative  for
$Out(\f)$ further by proving:

\begin{stoa} A solvable subgroup of $Out(\f)$ has a
finitely generated free abelian subgroup of index at
most $3^{5n^2}$.\end{stoa}

\noindent The rank of an abelian subgroup of $Out(\f)$ is $\leq
2n-3$ for
$n>1$ \cite{cv:moduli}.

\vskip.2cm

There is a reformulation of our Kolchin theorem  in terms of trees.
This is the
form in which we prove the theorem in this paper.

\begin{tree kolchin} 
For every finitely generated $UPG$ subgroup $\H$ of $Out(\f)$
there is a nontrivial simplicial
$\f$-tree with all edge stabilizers trivial that is fixed by all
elements of $\H$.\end{tree kolchin}

     Such a tree can be obtained from the marked filtered graph
produced by our Kolchin theorem by taking the universal cover and then
collapsing all edges except for the lifts of the highest edge $E_K$.
For a proof of the reverse implication, namely that Theorem 
\ref{kolchin for trees} implies Theorem \ref{graph kolchin} see Section
\ref{main proof}.

    We now outline the proof of Theorem \ref{kolchin for trees}.
The idea is
to find the common fixed tree using an iteration scheme. This
iteration takes place in the space $\x$ of very small simplicial
$\f$-trees (for a definition see Section \ref{very small trees}).
There is a natural (right) action of $Out(\f)$ on $\x$ (see
Section \ref{trees} for a review of the necessary background). In
the
first part of the paper we are primarily concerned with a study of
the
dynamics of the action of a $UPG$ automorphism $\oa$ on
$\x$. Specifically, we show in Theorem \ref{limitexists} that
under
iteration every tree $T\in\x$ converges to a tree
$T\oa^{\infty}\in\x$
(necessarily fixed
by the automorphism).
This fact is a consequence of Theorem
\ref{eventually polynomial}, which asserts that the sequence of
iterates of any $\gamma\in \f$ under a $UPG$ automorphism
$\oa$
eventually behaves like a polynomial (for a definition, see Section
\ref{polynomial sequences}; in particular the function $k\mapsto
length(\oa^k([\gamma]))$ coincides with a polynomial for large
$k$). In
addition, we show that the asymptotic behavior of the sequence
$\{\oa^k([\gamma])\}$ is largely determined by a finite number
of {\it
eigenrays} of $\oa$ that correspond to the eigendirections in the
linear
case.

Section \ref{kolchin} is the heart of the proof.
Let $T_0$ be a nontrivial simplicial $F_n$-tree with trivial edge
stabilizers such that the set of elliptic elements (a {\it free factor
system}) is $\H$-invariant and maximal among all $\H$-
invariant free
factor systems. For notational simplicity let us assume that $\H$
is
generated by two elements, $\oa_1$ and $\oa_2$.  Then consider
the
sequence $T_0,T_1,T_2,\cdots$ of simplicial trees defined
inductively
by $T_{i+1}=T_i\oa_1^\infty$ if $i$ is even and by
$T_{i+1}=T_i\oa_2^\infty$ if $i$ is odd. We then show that the
sequence is eventually constant and the tree thus obtained has the
desired properties. The first step consists of showing that a suffix
of $\oa_1$ or $\oa_2$ can be hyperbolic in at most one of the
trees in
the sequence. This claim is established by showing that, assuming
the
contrary, some element of $\H$ grows exponentially. The
argument is
reminiscent of the argument that the group generated by two Dehn
twists in intersecting curves contains an exponentially growing
mapping class. It follows from this first step that eventually all
suffixes of $\oa_1$ and $\oa_2$ are elliptic in $T_i$. The second
step
is to show that starting from this $T_i$ the set of elliptic elements
forms a
decreasing sequence. After establishing a chain bound on sets of
elliptic elements (see Proposition \ref{propC}), this implies that
the
set of elliptics in $T_j$ is independent of $j$ (for large $j$). By
the unipotent assumption, the vertex stabilizers of $T_j$ are fixed
by
$\H$ (rather than permuted) up to conjugacy (see Proposition
\ref{upgperiodic}). This however does not imply that $T_j$ is
fixed by $\oa_1$ and $\oa_2$ (for examples see Section
\ref{kolchin}). In the third step we examine the edge stabilizers of
$T_j$. If some are trivial and some nontrivial, then by collapsing
those with nontrivial stabilizer we obtain a tree that contradicts the
choice of $T_0$. If they are all nontrivial, we examine the term of
the sequence $\{ T_k\} $ that gave rise to such edges and again find a
larger proper free factor system invariant under $\H$. Therefore all
edges of $T_j$ have trivial stabilizer. In the fourth and final step
we observe that if $T_j$ is not fixed by $\oa_1$ and $\oa_2$, then in
the sequence $\{ T_k\}$ an equivalence class of edges gets short
compared to the average and there is again a larger proper free factor
system invariant under $\H$. This last step is reminiscent of the
development of Nielsen-Thurston theory where one finds invariant
curves of a mapping class by an iteration scheme in the moduli space
and picks out the curves that get short.

The key arguments in the paper focus not on discovering a ping-
pong
dynamics ($\H$ may well contain $F_2$) but on constructing an
element
in $\H$ of exponential growth. These are Proposition \ref{ugly},
Proposition
\ref{evnotgrow}, and Proposition \ref{nielsen}.

      After the breakthrough of E. Rips and the subsequent successful
applications of the theory by Z. Sela and others it became clear that
trees were the right tool for proving Theorem~\ref{graph kolchin}.
Surprisingly, under the assumption that $\H$ is finitely generated
(the case we are concerned with in this paper and that suffices for
the Tits Alternative), we only work with simplicial trees and the full
scale $\mathbb R$-tree theory is never used. However, its existence
gave us a firm belief that the project would succeed, and, indeed, the
first proof we found of the Tits Alternative used this theory. In a
sense, our proof can be viewed as a development of the program,
started by Culler-Vogtmann
\cite{cv:moduli}, to use spaces of trees to understand $Out(\f)$
in much
the same way that Teichm\" uller space and its compactification
were used by Thurston and others
to understand
mapping class groups.

The results of \cite{bfh:tits1} used
here are collected in Section \ref{unipotent section}. The reader interested
primarily in the arguments involving trees can read the present paper
independently of \cite{bfh:tits1}.

\section{${\f}$-Trees}\label{trees}

In this section, we collect the
facts about real $\f$-trees that we will need. This paper will only use
these facts for simplicial trees, but we record more general
results for later use.  

\subsection{Very small trees} \label{very small trees}

An $\f$-tree $T$ is {\it very small} \cite{cl:verysmall} if it is
minimal (i.e. it does not have any proper invariant subtrees),
nondegenerate (i.e. it is not a point), all edge stabilizers are
trivial or primitive cyclic, and for each $1\neq\gamma\in \f$ the
subset $Fix_T(\gamma)$ of $T$ fixed by $\gamma$ is either empty, a
point, or an arc. The set of all projective classes of very small
$\f$-trees is denoted by $\xx$ and the subset of $\xx$ consisting of
the projective classes of simplicial trees is denoted by $\x$. Both are
topologized via the embedding $\theta:\xx\to {\mathbb P}^{\cal C}$ into
the infinite-dimensional projective space, where $\Cal C$ is the set
of all conjugacy classes in $\f$ and $\theta (T):[\gamma]\mapsto
\ell_T(\gamma)$ ($\ell_T(\gamma)$ is the translation length of
$\gamma$ in $T$). See \cite{cm:trees} for a proof that $\theta$ is
injective.

The automorphism group $Aut(\f)$ acts naturally on $\xx$ on the right
by changing the marking. In terms of the length functions, the action
is given by $\ell_{T\oa}([\gamma])=\ell_T(\oa([\gamma]))$. Inner
automorphisms act trivially and we have an action of
$Out(\f)=Aut(\f)/Inn(\f)$. There is a natural invariant decomposition
of $\x$ into open simplices. The space $\xx$ can be identified
\cite{cl:verysmall} \cite{cm:trees} with Culler-Morgan's
compactification \cite{cm:trees} of Culler-Vogtmann's Outer Space
\cite{cv:moduli}.

\subsection{Bounded cancellation constants}\label{bccsec}

We will often need to compare the length of the same element
of $\f$ in different
$\f$-trees. The existence of bounded cancellation constants will
usually suffice for this job.

\begin{defn} The {\it bounded cancellation constant} 
of an $\f$-map $f:T'\to T$,
denoted $BCC(f)$, is the least upper bound of numbers $B$ with the
property that there exist points $a,b,c\in T'$ with $b$ on the segment
$[a,c]$ so that the distance between $f(b)$ and the segment
$[f(a),f(c)]$ is $B$.\end{defn}

In \cite{co:bcc} Cooper showed that if both $T$ and $T'$ are free
simplicial and minimal and $f$ is $PL$, then $BCC(f)$ is
finite. The bound given by Cooper depends on the Lipschitz constants
of $f$ and of an $\f$-map $T'\to T$.

Below we generalize Cooper's result to the case that the target tree $T$ is
very small. For a map $f$ between metric spaces we denote by $L(f)$ the
Lipschitz constant of $f$, i.e.
$$L(f):=sup\{ {{d_{T'}(f(a),f(b))}/ {d_T(a,b)}}|(a,b)\in T\times
T, a\not= b\}.$$

\begin{lem} \label{bcccomp} 
Suppose $f:T''\to T'$ and $g:T'\to T$ are $\f$-maps
between minimal $\f$-trees. Then
\begin{enumerate}
\item 
$BCC(g)\le BCC(gf)$
\item 
$BCC(gf)\le BCC(g)+ L(g)BCC(f)$
\end{enumerate}
\end{lem}

\begin{proof} 
(1) follows directly from definition.

(2)
Choose $a,b,c\in T''$ with
$b\in [a,c]$, and let $b'$ be the point in $[f(a),f(c)]$ closest to
$f(b)$. Then 
\begin{equation*}
\begin{split}
d(gf(b),[gf(a),gf(c)])&\allowbreak\le
d(gf(b),g(b'))+\allowbreak d(g(b'),[gf(a),gf(c)])\\ &\le 
L(g)BCC(f)+BCC(g).\end{split}\end{equation*}\end{proof}

\begin{defn} The {\it covolume} of a free simplicial and minimal $\f$-tree $T$,
denoted $cov(T)$, is the sum of the lengths of edges in $T/\f$.\end{defn}

\begin{prop} \label{splittech} Suppose $f:T_1\to T$ is an $\f$-map
between free simplicial and minimal $\f$-trees $T_1$ and $T$.
Assume $f$ is linear on each edge of $T_1$. Then
$$
BCC(f)\le L(f)cov(T_1)-cov(T)
$$
\end{prop}

\begin{proof} Represent $f$ as the composition $hgf_k\cdots f_1$ where each
$f_i:T_i\to T_{i+1}$ is a fold (see \cite{st:folding}), $g$
collapses some orbits of edges, and $h$ is a homeomorphism, linear on
each edge, with $L(h)=L(f)$ (see \cite[page 452]{bf:bounding}). 
Note that $L(g)=L(f_i) =1$,
$BCC(h)=BCC(g)=0$, and $BCC(f_i)=cov(T_i)-cov(T_{i+1})$. Then we use
Lemma
\ref{bcccomp}.
\begin{equation*}
\begin{split}
BCC(f)&=BCC(hgf_k\cdots f_1)\\ &\leq BCC(h)+L(h)BCC(gf_k\cdots f_1)\\ &\leq
BCC(h)+L(f)(BCC(g)+BCC(f_k)+\cdots+BCC(f_1))\\ &=L(f)(BCC(f_k)+\cdots+BCC(f_1))
\\ &=
L(f)(cov(T_1)-cov(T_{k+1}))\\ &\leq
L(f)cov(T_1)-cov(T)\end{split}\end{equation*}
\end{proof}

\begin{prop} \label{finite bcc}
Suppose $f:T_1\to T$ is a Lipschitz $\f$-map
from a free simplicial $\f$-tree $T_1$ to a very small simplicial
$\f$-tree $T$. Assume $f$ is linear on each edge of $T_1$. Then
$BCC(f)<\infty$.
\end{prop}

\begin{proof} Represent $f$ as a composition of folds and apply Lemma
\ref{bcccomp}.\end{proof}

The following generalization is not used in this paper, but will be in
\cite{bfh:tits3}.

\begin{prop} For every tree $T\in\xx$, every very small simplicial
$\f$-tree $T_1$ any $\f$-map $T_1\to T$ that is linear on edges has
finite $BCC$.
\end{prop}

\begin{proof} By Lemma \ref{bcccomp} 
it suffices to prove the proposition in the case that $T_1$ is the
universal cover of a rose. Further, it suffices to
construct an $F_n$-map $T_1\to T$ with finite $BCC$ (since any two
such maps are within bounded distance from each other). There is an
embedding $\phi$ of $\x$ into the space of very small trees. Indeed,
let $\{ x_1,\cdots,x_{\n}\}$ be a basis for $\f$ and let $\cal P$
denote the elements of word length at most 2. If $T$ is a nontrivial
$\f$-tree then the lengths of the elements of $\cal P$ can't all be 0
\cite{cv:moduli}. Thus, the set of all very small trees $T$ such that
the sum of the lengths of the elements of $\cal P$ equals 1 is
homeomorphic to $\x$. Let $T_1$ be the universal cover of a rose 
in $\phi(\x)$. There is a continuous choice of base
point for each $T\in\phi(\x)$
\cite{rs:deformations},\cite{tw:fixed}. Let $f_T:T_1\to T$ be the map
that takes the vertex of $T_1$ to the base point of $T$ and is linear
on the edges of $T_1$. Since the topology on $\phi(\x)$ is the same as
the based length function topology \cite{ab:lengthfunctions}, the Lipschitz
constant $L(f_T)$ varies continuously. Since the topology on
$\phi(\x)$ is the same as the equivariant Gromov-Hausdorff topology,
$BCC$ is lower semi-continuous \cite{fp:trees}. By Proposition
\ref{splittech}, $BCC(f_T)\le L(f_T)cov(T_1)$ if $T$ is minimal free
simplicial. So, the proof now follows from the above observations
together with the fact that every very small tree is the limit of free
simplicial and minimal trees \cite[Theorem
2.2]{bf:outerlimits}.\end{proof}

\subsection{Free factor systems}\label{free factor systems}


Our next goal is to prove that chains of the sets of elliptic elements
in very small $F_n$-trees are bounded. To develop notation we first handle
the special case of free factor systems, which corresponds to restricting
to simplicial trees with trivial edge stabilizers. There is some overlap
between this section and \cite{bfh:tits1}.

Let $\Cal N$ denote the set of finite nonincreasing sequences in $\mathbb
N$. We allow the empty sequence. Well order $\Cal N$
lexicographically. For example, $5,3,3,1 > 4,4,4,4,4,4 > 4>\emptyset$.
In the cases that we consider, the sum of the elements in the set will
be no more than $\n$.  Thus, the sequence $\n$ will be the largest
element that we will consider and $\emptyset$ the smallest.

\begin{defn} If $F$ is a subgroup of ${\f}$ then let $[F]$
denote the set of all subgroups of ${\f}$ that are conjugate to $F$,
i.e. the conjugacy class of
$F$. A set
$\Cal F$ of free factors of ${\f}$ is a {\it free
factor system} if there is a free factor of ${\f}$ of the form
$F_1*...*F_k$ such that ${\Cal F} =
[F_1]\cup[F_2]\cup\cdots\cup[F_k]$. For convenience, we
will always require that $<1>\in\Cal F$. Equivalently, a free factor
system is the set of point stabilizers of a simplicial ${\f}$-tree with
trivial edge stabilizers.  The {\it complexity} of
$\Cal F$ is the element of $\Cal N$ obtained by arranging the
positive numbers among
$rank(F_1), ...  , rank(F_k)$ in nonincreasing order.
$\Cal F$ is {\it proper} if its complexity is less than $\n\in\Cal N$.
\end{defn}

\begin{lem} If $\Cal F$ and $\Cal F'$ are two free factor
systems then $\{ F\cap {F'} | F\in {\Cal F}, F'\in {{\Cal F}'}\}$ is a free
factor system.  \end{lem}

Denote this free factor system by ${\Cal F} \wedge {\Cal F}'$.

\begin{proof} Let $T_{\Cal F}$ denote a simplicial tree with trivial
edge stabilizers and vertex stabilizers $\Cal F$. Let $T_{\Cal F'}$
denote a similar tree with respect to $\Cal F'$. For $F\in \Cal F$,
consider the action of $F$ on $T_{\Cal F'}$.  This gives a simplicial
$F$-tree with vertex groups $\{ F\cap {F'}^c | F'\in {\Cal F'}, c\in
{\f}\}$ and trivial edge groups. Use this tree to blow up
\cite{rj:blowup} the orbit of the vertex of $T_{\Cal F}$ stabilized by
$F$.  We obtain a simplicial ${\f}$-tree with trivial edge stabilizers
and vertex stabilizers ${\Cal F \wedge \Cal F'}$.\end{proof}

\begin{notn} If $\Cal H$ and $\Cal H'$ are two sets of
subsets of ${\f}$ then we write $\Cal H\preceq\Cal H'$ if each $H\in\Cal
H$ is contained in some
$H'\in\Cal H'$. If also $\Cal H\not=\Cal H'$ then we
write
$\Cal H\prec\Cal H'$.\end{notn}

\begin{prop} \label{ffs chain bound} Let $\Cal F$ and $\Cal F'$ be two
free factor systems. If $\Cal F\preceq\Cal F'$ then $$Complexity({\Cal
F})\le Complexity({\Cal F'}).$$ If additionally $\cup{\Cal F}\not=
\cup{\Cal F'}$ then
${\Cal F}\prec\Cal F'$ and $$Complexity({\Cal F})<Complexity({\Cal
F'}).$$\end{prop}

\begin{proof} If $F$ and $F'$ are free factors on ${\f}$ such that
$F\subseteq F'$ then $rank(F)\le rank(F')$ with equality if and only if
$F=F'$. The lemma now follows easily.\end{proof}

\begin{cor} $$Complexity({\Cal F}\wedge {\Cal F'})\le
Min\{ Complexity({\Cal F}),Complexity({\Cal F'})\}.$$\end{cor}

\begin{lem} \label{min} Let $\Cal H$ be a (possibly infinite) set of
subsets of ${\f}$.  Then there is a free factor system $\Cal F({\Cal H})$ of minimal complexity such that $\Cal H \preceq
\Cal F({\Cal H})$.  Further, this system is unique.\end{lem}

\begin{proof} Clearly there is such a system, call it $\Cal F$. If
$\Cal F'\not= \Cal F$ is another then so is
$\Cal F \wedge \Cal F'$ but of smaller complexity. \end{proof}

\begin{cor} $$Complexity\big({\Cal F}({{\Cal H} \cup {\Cal H'}})\big)\ge
Max\{ Complexity\big({\Cal F}({\Cal H})\big),Complexity\big({\Cal F} ({\Cal
H'})\big)\}.$$\end{cor}

\begin{notn} \label{closure}
Let $\partial {\f}$ denote the Hopf boundary
\cite{hopf:ends} of ${\f}$ (which agrees with the Gromov boundary in this
case). If ${\f}$ is represented as the fundamental group of a graph $\G$, then
$\partial {\f}$ may be identified with the space of geodesic rays in the
universal cover $\tilde \G$ of $\G$ where we identify two rays
if they eventually coincide. For a finitely generated subgroup
$F\subseteq{\f}$, inclusion is a quasiisometric embedding (see Lemma
\ref{qiembedded} below) and so we may identify
$\partial F$ with a subset of $\partial {\f}$. If $F$ is represented as
subgraph $\Delta$ of $\G$ then $\partial F$ may be identified with
the subspace of geodesic rays that are eventually in the preimage of
$\Delta$ in $\tilde\G$. Let $\overline F$ denote the subset
$F\cup\partial F$ of ${\f}\cup\partial{\f}$. For a set $\Cal F$ of
finitely generated subgroups of ${\f}$ let $\overline{\Cal F}$ denote $\{
\overline F | F\in\Cal F\}$. If $\Cal H$ and $\Cal H'$ are two sets of
subsets of ${\f}\cup\partial{\f}$ then we write $\Cal H\preceq\Cal H$ if each
$H\in\Cal H$ is contained in some $H'\in\Cal H'$. If also $\Cal
H\not=\Cal H'$ then we write $\Cal H\prec\Cal H'$.\end{notn}

For a proof of the following lemma, in far greater generality, see
\cite{short:quasiconvex}.

\begin{lem} \label{short}
Let $H, H'$ be finitely generated subgroups of ${\f}$. 
\begin{enumerate}
\item[(1)] \label{qiembedded} The inclusion
$H\subseteq \f$ is a quasiisometric embedding.

\item[(2)] \label{capb}
$\bd
F\cap \bd F' = \bd(F
\cap F')$.\qed\end{enumerate}\end{lem}

Using Lemma \ref{short}, a proof similar to that of Lemma
\ref{min} establishes:

\begin{lem} Let $\Cal H$ be a set of subsets of ${\f}
\cup
\partial {\f}$. Then there is a unique free factor system
$\Cal F ( {\Cal H})$ of minimal
complexity such that $\Cal H\preceq\overline{\Cal F(\Cal
H)}$.\end{lem}


\subsection{A chain bound for vertex systems}

Since some of our arguments will proceed by restricting outer
automorphisms to point stabilizers of trees in $\xx$, it is important
to get a precise picture of these stabilizers. In the case where $T$
is simplicial with trivial edge stabilizers, the set of point
stabilizers is a free factor system, and we have analyzed these in
Section \ref{free factor systems}.

\begin{defn}
A {\it vertex group} is a point stabilizer of a tree in $\xx$.
For an $\f$-tree $T$, $\ve(T)$ denotes the collection of its point
stabilizers, and $\cup\ve(T)$ is the set of elliptic group elements.\end{defn}

In this section, we show the existence of a bound for the length of
sequence of inclusions of vertex groups (Proposition
\ref{vertexchainbound}) or more generally vertex systems (Proposition
\ref{propC}). 

In the case of simplicial trees, the following theorem is
established by an easy Euler characteristic argument. The following
generalization to $\mathbb R$-trees due to Gaboriau and Levitt uses
more sophisticated techniques.

\begin{thm}\cite{gl:rank}\label{gl} 
Let $T\in\xx$. There is a bound depending only on $\n$ to the number
of conjugacy classes of point and arc stabilizers. The rank of a point
stabilizer is no more than $\n$ with equality if and only if $T/\f$ is
a rose and each edge of $T$ has infinite cyclic
stabilizer.\qed\end{thm}

\begin{prop} \label{vertexchainbound} There is a bound (depending only
on $\n$) to the length of any chain of proper inclusions of
vertex groups.\end{prop}

\begin{proof} Let $V \supset V' \supset V''$ be a chain of proper
containments of vertex groups with corresponding trees $T$, $T'$, and
$T''$. We will show that either $rank(V) > rank(V'')$ or $\n \ge
rank(H_1(V/\ll V''\gg )) > rank(H_1(V/\ll V'\gg ))$.

Let $T'_{V}$ and $T''_{V}$ be minimal $V$-subtrees of
$T'$ and $T''$ respectively. Since the vertex groups of $T'_V$ are
precisely the intersection of the vertex groups of $T'$ with $V$, we
see that $T'_V$ has a vertex labelled $V'$ and so $rank(V) \ge
rank(V')$. Similarly, $rank(V')\ge rank(V'')$. If $rank(V) >
rank(V'')$ then we are done, so assume these ranks are equal.

Using Theorem \ref{gl}, the only remaining possibility is that
orbit spaces of both $T'_V$ and
$T''_V$ are roses of circles with all edges labeled by infinite cyclic
groups. The number $k(T'_V)$ of orbits edges of $T'_V$ may be computed as
$rank(H_1(V/\ll V'\gg )$. (Indeed, in general if $S$ is a $V$-tree, $X$
its orbit space
$S/V$, and $K=Kernel(V\to
\pi_1(X))$ then $K=<\cup\ve(S) >$.)
Since there is an epimorphism
$$V/\ll V''\gg \longrightarrow V/\ll V'\gg $$ we have that $\n \ge k(T''_V) \ge
k(T'_V)$. We will show that if $k(T''_V) = k(T'_V)$ then $V' = V''$, a
contradiction.

Consider the morphism $\phi:T''_V\to T'_V$ that sends the vertex $v''$
labelled $V''$ to the vertex $v'$ labelled $V'$.  This map is well
defined for if $e''$ is the edge from $v''$ to $gv''$ with stabilizer
$E''$, then $E''=V''\cap {V''}^g \subseteq V'\cap {V'}^g
=Stabilizer(\phi e'')$. Thus, $T'_V$ is obtained from $T''_V$ by a
finite number of folds (after perhaps first subdividing $T''_V$)
\cite[page 455]{bf:bounding}.  An inspection of the types of folds 
\cite[pages 452--3]{bf:bounding} reveals that, in this situation, a sequence of folds 
cannot change the vertex groups without decreasing the first Betti number
of the quotient graph or increasing the rank of an edge
stabilizer. 
\end{proof}

\begin{com} Notice that a vertex group of a vertex group is not necessarily
a vertex group. This is apparent from Proposition \ref{vertexchainbound} and
the fact that there are $\f$-trees with vertex groups of rank that of
$\f$.\end{com}

\begin{lem} \label{maxisstab} Let $T\in\xx$. 
Then $\cup\ve(T)$ is the union of the maximal groups in $\ve(T)$ each
of which is a point stabilizer in $T$.\end{lem}

\begin{proof} The lemma will follow if we show that every
group
$V\subseteq\cup\ve(T)$ fixes a point in $T$. If $V$ is finitely
generated then, by \cite[page 65]{se:trees}, the restriction of the action of
$V$ to $T$ can have a trivial length function only if there is a
global fixed point. Thus we may assume that $V$ is not finitely
generated. Hence,
$V$ is an increasing union of noncyclic finitely generated subgroups
each contained in some point stabilizer. A noncyclic group fixes at
most one point of $T$ hence $V$ fixes a point.\end{proof}

\begin{prop} \label{propC} There is a bound depending only on $\n$ to the
length of a sequence of inclusions
$$\cup\ve(T_0)\subsetneq\cup\ve(T_1)\subsetneq\cdots\subsetneq
\cup\ve(T_N)$$ where each $T_i\in\xx$.\end{prop}

\begin{proof} Let $l$ be a bound on the length of a chain of proper inclusions
of vertex groups of trees in $\xx$. The existence of $l$ is
guaranteed by Proposition \ref{vertexchainbound}. Let ${\Cal M}_i$ denote the
set of conjugacy classes of maximal groups in $\cup\ve(T_i)$. By Lemma
\ref{maxisstab}, ${\Cal M}_i$ consists of conjugacy classes of point stabilizers.
For $M\in{\Cal M}_i$, let $\cup M\subseteq \f$ denote the set of elements
represented by $M$. Let ${\Cal A}_i$ denote the subset of the power
set $\Cal P(\f)$ of $\f$ given by
$\{\cup M|M\in {\Cal M}_i\}$. Let $k$ be a bound to $|{\Cal A}_i|$
(see Theorem \ref{gl}). We will show
that
$N<l^k$.

\begin{sublem} \label{comb} For every pair of integers $k,l>0$ the following
holds with
$n=l^k$. Let ${\Cal A}_0,{\Cal A}_1,\dots,{\Cal A}_n$ be subsets of the power set
${\Cal P}(X)$ of a fixed set $X$. Assume that
\begin{itemize}
\item if $A,A'\in{\Cal A}_i$ with $A\subseteq A'$ then $A=A'$,
\item $|{\Cal A}_i|\leq k$ for all $i$, and 
\item ${\Cal A}_i\preceq {\Cal A}_{i+1}$ for all $i$.
\end{itemize}
Then one of the following two possibilities occurs.
\begin{enumerate}
\item There are $A_i\in {\Cal A}_i$, $i=1,2,\dots,n$, such that
$$A_0\subseteq A_1\subseteq\dots\subseteq A_n$$
and at least $l$ of these inclusions are proper.
\item For some $i$, ${\Cal A}_i={\Cal A}_{i+1}$.
\end{enumerate}
\end{sublem}

\begin{proof} Induction on $k$. The case $k=1$ is clear. Now suppose that
the lemma is true for $k-1$. Choose arbitrary $A_i\in {\Cal A}_i$ with
$A_0\subseteq A_1\subseteq \dots \subseteq A_n$. Consider $l$ chains
of inclusions of length $l^{k-1}$
$$A_{jl^{k-1}}\subseteq A_{jl^{k-1}+1}\subseteq \dots\subseteq
A_{(j+1)l^{k-1}},\ j=0,1,2,\dots,l-1.$$ If each chain contains a
proper inclusion, then (1) holds. If not, then one of these chains,
say the one with $j=0$, consists of equalities. Now, remove $A_i$ from
${\Cal A}_i$, $i=1,2,\dots,l^{k-1}$. By the first bullet, the new collection
satisfies the inductive hypothesis.\end{proof}

We now continue the proof of Proposition \ref{propC}. Assume $N=l^k$.
The hypotheses of Sublemma
\ref{comb} are satisfied. According to the Sublemma,
there is $i$ so that ${\Cal A}_i={\Cal A}_{i+1}$, since (1) is
impossible by our choice of $l$. But now we have
$\cup\ve(T_i)=\cup\ve(T_{i+1})$ since $\cup\ve(T_i)=\cup{\Cal A}_i$.
\end{proof}


\section{Unipotent polynomially growing outer automorphisms}\label{unipotent}

We now bring outer automorphisms into the picture. We will
consider a class of outer automorphisms that is analogous to the class of
unipotent matrices.  First we review the linear algebra of unipotent
matrices.

\subsection{Unipotent linear maps}

\begin{prop-def}
Let $R=\mathbb Z$ or $\mathbb C$,
and let $V$ be a free $R$-module of finite rank. We say that an
$R$-module endomorphism $F:V\to V$ is {\it unipotent} if the following
equivalent conditions are satisfied:
\begin{enumerate}
\item $V$ has a
basis with respect to which $F$ is upper triangular with 1's on the diagonal.
\item $(Id-F)^n=0$ for some $n>0$.
\end{enumerate}
\end{prop-def}

\begin{proof} It is clear that (1) implies (2). To see that (2) implies
(1), assume that $(Id-F)^n=0$. We may assume that $W:=Im(Id-F)^{n-1}\neq 0$.
The restriction of $Id-F$ to the submodule $W$ is 0, and hence each $0\neq w\in
W$ is fixed by $F$. In the case $R=\mathbb Z$ pass to a primitive submultiple if
necessary to conclude that $V$ always contains an $F$-fixed basis element $v$.
The proof now concludes by induction on $rank(V)$ using the observation that
the induced homomorphism $F':V/<v>\to V/<v>$ also satisfies $(Id-F')^{n}=0$.
\end{proof}

\begin{cor} \label{uprestrict} Let $R=\mathbb Z$ or $\mathbb C$. Let
$F:V\to V$ be an $R$-module endomorphism, and let $W$ be an $F$-invariant
submodule of $V$ which is a direct-summand of $V$. Then $F$ is unipotent if
and only if both the restriction of $F$ to $W$ and the induced endomorphism on
$V/W$ are unipotent.\end{cor}

\begin{proof} Evident, if we use (2) in $\Longrightarrow$ and (1) in
$\Longleftarrow$.\end{proof}

\begin{cor}\label{periodicisfixed} 
Let $F:V\to
V$ be unipotent. If $x\in V$ is $F$-periodic, i.e. if $F^m(x)=x$ for some
$m>0$, then $x$ is $F$-fixed, i.e. $F(x)=x$.\end{cor}

\begin{proof} First assume that $R=\mathbb C$. We may assume that
$V=span(x,F(x),\cdots ,\allowbreak F^{m-1}(x))$. Let $e_1,e_2,\dots,e_m$
be the standard basis for $\mathbb C^m$. There is a surjective linear map
$\pi:\mathbb C^m\to V$ given by $\pi(e_i)=F^{i-1}(x)$, and $F$ lifts to
the linear map $\overline F:\mathbb C^m\to \mathbb C^m$, $\overline
F(e_i)=e_{i+1\mod m}$. For $\lambda\in\mathbb C$, the generalized
$\lambda$-eigenspace is defined to be
$$
\{ x\in \mathbb C^m| (\lambda I-\overline F)^m(x)=0\}.
$$
The linear map
$\pi$ must map the generalized
$1$-eigenspace onto $V$ (and all other generalized eigenspaces to 0). Since
this space is one-dimensional (and equals the $1$-eigenspace of $\overline
F$), it follows that $dim(V)\leq 1$ and $F(x)=x$.

If $R=\mathbb Z$, just tensor with $\mathbb C$.
\end{proof}

\begin{cor}\label{summands} 
Let $F:V\to V$ be unipotent. If $W$ is a direct
summand which is periodic (i.e. $F^m(W)=W$ for some $m>0$), then $W$
is invariant (i.e. $F(W)=W$).\end{cor}

\begin{proof} The restriction of $F^m$ to $W$ is unipotent, so there is a basis
element $x\in W$ fixed by $F^m$. By Corollary 3.3, $F(x)=x$. The proof
concludes by induction on $rank(W)$.\end{proof}

\begin{prop} \label{unichar}
Let $F\in GL_n(\mathbb Z)$ have all eigenvalues on the unit circle (i.e.
$F$ grows polynomially). If the image of $F$ in $GL_n(\mathbb Z/3)$
is trivial, then $F$ is unipotent.\end{prop}

\begin{proof}
We first argue that some power $A^N$ of $A$ is unipotent, i.e. that
all eigenvalues of $A$ are roots of unity. Choose $N$
so that all eigenvalues of $A^N$ are close to 1. Then $tr(A^N)$ is an
integer close to $n$, and thus all eigenvalues of $A^N$ are equal to 1.

Let $f=f_1^{n_1}\cdots f_m^{n_m}$ be the minimal polynomial for $A$
factored into irreducibles in $\mathbb Z[x]$. Let $A_i=f_i^{n_i}(A)$
and $K_i=Ker(A_i)$.  First note that each $K_i\neq 0$. For example,
$Im(A_2 A_3\cdots A_m)\subset K_1$ but $A_2A_3\cdots A_m\neq 0$ since
$f$ is minimal. If $A$ is not unipotent, then some $f_i$, say $f_1$, is
not $x-1$. Thus $f_1$ is the minimal polynomial for a nontrivial root
of unity and so it divides $1+x+x^2+\cdots +x^{r-1}$ for some
$r>1$. The matrix $I+A+A^2+\cdots +A^{r-1}$ has nontrivial kernel
(since its $n_1^{st}$ power vanishes on $K_1$). It follows that there is
a nonzero integral vector $v$ such that $A^r(v)=v$ but $A(v)\neq
v$. Then $Fix(A^r)$ is a nontrivial direct summand of $\mathbb Z^n$,
the restriction of $A$ to this summand is nontrivial and periodic, and
the induced endomorphism of $Fix(A^r)\otimes \mathbb Z/3$ is identity.
This contradicts the standard fact that the kernel of $GL_k(\mathbb Z)\to
GL_k(\mathbb Z/3)$ is torsion-free.
\end{proof}

\subsection{Relative train tracks}

Techniques of this paper strongly depend on finding good
representatives for {\it polynomially growing} outer automorphisms.

\begin{defn} An outer automorphism $\oa\in Out(\f)$ is $PG(\f)$ (or just $PG$)
if for each conjugacy class $[\gamma]$ in $\f$ the sequence
of (reduced) word lengths of $\oa^i([\gamma])$ is bounded above by a
polynomial. \end{defn}

We start by recalling the representatives for $PG$ automorphisms found
in \cite{bh:tracks}.

\begin{thm} \label{tracks} \cite{bh:tracks} 
Every $PG$ automorphism $\oa\in Out({\f})$ has a
representative as a homotopy equivalence $\tmm:\ttt\to\ttt$ on a marked graph
$\ttt$ such that \begin{enumerate}
\item the map $\tmm$ sends vertices to vertices and edges to immersed
nontrivial edge paths.
\item There is a filtration 
$\emptyset=G_0\subsetneq G_1\subseteq\cdots\subsetneq G_K=G$ of $G$ by
$\tmm$-invariant subgraphs such that for every edge $E\in
\overline{G_i\setminus G_{i-1}}$ 
the edge path $\tmm(E)$ crosses exactly one edge in $\overline{G_i\setminus
G_{i-1}}$ and it crosses that edge exactly once.
\item If $\Cal F$ is an $\oa$-invariant free factor system, we can
arrange that $\Cal F$ is represented by some $G_r$. If $\oa$ is the identity
on each conjugacy class in $\Cal F$, we can arrange that $\tmm=Id$ on
$G_r$.
\end{enumerate}
\end{thm}

\definition The representative in Theorem \ref{tracks} is called a
{\it relative train track ($RTT$) representative} for $f$. 
\enddefinition

\begin{notn}
All paths in graphs and trees will have endpoints in the vertex set. If
$\gamma$ is a path, $[\gamma]$ will denote the unique immersed path
homotopic to $\gamma$ rel endpoints. When the endpoints of $\gamma$
coincide, we say that $\gamma$ is a {\it based loop}. When $\gamma$ is
an (unbased) essential loop, $[\gamma]$ will denote the unique immersed
loop freely homotopic to $\gamma$. We make standard identifications
between homotopy classes of based loops with elements of the
fundamental group and between homotopy classes of loops and conjugacy
classes in the fundamental group. 
\end{notn}

\subsection{Unipotent representatives} \label{unipotent section}

We now introduce $UPG$ automorphisms - the objects of central
importance in this paper.

\definition An outer automorphism 
is a $UPG(\f)$ (or just $UPG$) automorphism if it is $PG(\f)$ and its
action in $H_1(\f;\mathbb Z)$ is unipotent.\enddefinition

We now recall a special case of an improvement of $RTT$ representatives from
\cite{bfh:tits1}.

Recall that if $f:G\to G$ is a
RTT representative and $z$ is an edge path in $G$ then we write $[z]$
for the geodesic homotopic rel endpoints to $z$. We also write $z=x\cdot y$
for edge paths $x$ and $y$ in $G$ if $[f^n(z)] = [f^n(x)][f^n(y)]$ for all
$n\ge 0$ and say that $z$ ``splits''.

\begin{defn} Let $f:G\to G$ be an $RTT$ representative. 
A path $\tau$ in $G$ with endpoints in the vertex set is {\it Nielsen} if 
$[f(\tau)]=[\tau]$. An
{\it exceptional path} in $G$ is a path of the form

$\bullet$ $E_i\tau^m E_i^{-1}$ provided $\tau$ is a nontrivial Nielsen path
and 
$f(E_i)=E_i\tau^p$ for some $m,p\in\mathbb Z$, $m\neq 0$, or

$\bullet$ $E_i\tau^m E_j^{-1}$ provided $\tau$ is a Nielsen path,
$i\neq j$, $f(E_i)=E_i\tau^p$, and $f(E_j)=E_j\tau^q$ for some
$m,p,q\in\mathbb Z$.\end{defn}

The following theorem follows easily from Theorem 6.8*** and Lemma 6.24***
of \cite{bfh:tits1}.

\begin{thm}(\cite{bfh:tits1})\label{improved} 
Suppose that
${{\cal O}} \in Out(F_n)$ is a $UPG$-automorphism, that
${{\cal F}}$ is an ${{\cal O}}$-invariant \ffs.
Then there is an $RTT$ representative $f:G\to G$
of
${{\cal O}}$ with the following properties.
\begin{enumerate}
\item ${{\cal F}} = {{\cal F}}(G_r)$ for some filtration
element $G_r$.
\item Each $G_i$ is the union of $G_{i-1}$ and a single edge
$E_i$ satisfying $f(E_i) = E_i \cdot u_i$ for some closed path $u_i$
that crosses only edges in $G_{i-1}$.
\item If $\sigma$ is any path with endpoints at
vertices, then there exists $M
= M(\sigma)$ so that for each $ m \ge M$,
$[f^m(\sigma)]$ splits into
subpaths that are either single edges or exceptional
subpaths.
\item $M(\sigma)$ is a bounded multiple of the
edge length of $\sigma$.
\item There is a uniform constant $C$ so that if $\omega$ is a
closed path in $G$ that is not a Nielsen path and
$\sigma=\alpha\omega^k\beta$ is an immersed path, then at most $C$
copies of $[f^m(\omega)]$ are cancelled when
$[f^m(\alpha)][f^m(\omega^k)][f^m(\beta)]$ is tightened to
$[f^m(\sigma)]$.
\end{enumerate}
\end{thm}

\begin{defn} An $RTT$ representative $f$ satisfying 1-5 above is a
{\it unipotent representative} or a $UR$. The based loops $u_i$ are
{\it suffixes} of $f$.\end{defn}

Note that (2) can be restated as $$[f^k(E_i)]=E_iu_i[f(u_i)]\cdots
[f^{k-1}(u_i)]$$ for all $k>0$. The immersed infinite ray
$$R_i=E_iu_i[f(u_i)]\cdots
[f^{k-1}(u_i)]\cdots$$
is the {\it eigenray} associated to $E_i$. Lifts of $R_i$ to the universal
cover of $G$ are also called eigenrays. The subpaths $[f^m(u_i)]$ of $R_i$
are sometimes referred to as {\it blocks}.

For example, the map $f:G\to G$ on the rose with two petals labelled $a$
and $b$ given by $f(a)=a$, $f(b)=ba$ is a $UR$. For $\omega=ba^{-10}bab^{-1}$
we may take $M(\omega)=10$ in (3), since $[f^{10}(\omega)]=b\cdot (bab^{-1})$
is a splitting into an edge and an exceptional (Nielsen) path. The map
given by $a\mapsto a$, $b\mapsto ba$, $c\mapsto cba^{-1}$ on the rose
with three petals is not a $UR$ since $\omega=cba^{-1}$ does not
eventually split as in (3). Replacing $ba^{-1}$ by $b'$ yields
a $UR$ of the same outer automorphism.

\begin{defn} Let $f:G\to G$ be a $UR$.
The {\it height}
of an edge-path in $G$ is the smallest $m$ such that the path is contained in
$G_m$. A {\it topmost edge} in an edge-path of height $m$ is an occurrence
of $E_m$ or $E_m^{-1}$ in the edge-path.\end{defn}

Many arguments are inductions on height. The inductive step is based on
the observation that a path of height $m$ splits at the initial (terminal)
endpoints
of each occurrence of $E_m$ ($E_m^{-1}$).

\section{The dynamics of $UPG$ automorphisms}

In this section we examine the dynamics of the action of $UPG$
automorphisms on conjugacy classes, free factor systems, and the space
$\xx$ of very small $\f$-trees.

\subsection{Polynomial sequences}\label{polynomial sequences}

In this section we show that the sequence of iterates of a path under
a $UPG$ automorphism behaves like a polynomial and use this to prove
Theorem \ref{limitexists} which is fundamental in our approach.

\begin{defn}
Let $G$ be a graph. A sequence $\{ A_k\}_{k=k_0}^{\infty}$ of
immersed paths in $G$ is said to be a {\it polynomial sequence} if
it can be obtained from constant sequences of paths by finitely many
operations described below. \begin{enumerate}
\item ({\it reindexing and
truncation}): $A_k=B_{k+k'}$ for a polynomial sequence $\{
B_k\}_{k=k_1}^{\infty}$ for some $k_1\leq k_0+k'$,
\item ({\it inversion}): $A_k$ is the inverse of $B_k$, and $\{
B_k\}_{k=k_0}^{\infty}$ is a polynomial sequence, \item ({\it
concatenation}): $A_k=B_k C_k$, where $\{ B_k\}$ and $\{ C_k\}$ are
polynomial sequences, (and no cancellation occurs in $B_kC_k$), and
\item ({\it integration}): $A_k=B_{k_0} B_{k_0+1}\cdots B_k$, where
$\{ B_k\}$ is a polynomial sequence (and again no cancellation
occurs).
\end{enumerate}\end{defn}

For example, in the standard
rose, sequences
$\{ AB^kC\}$ and $\{ ABAB^2AB^3\allowbreak\cdots AB^k\}$ are polynomial.

\begin{thm} \label{eventually polynomial}
Let $f:G\to G$ be a $UR$ representative of a $UPG$ automorphism. Let $\sigma$
be a path in $G$ with endpoints in the vertex set of $G$. Then there
is $k_0>0$ such that the sequence $\{ [f^k(P)]\}_{i=i_0}^{\infty}$ is
polynomial.  \end{thm}

\begin{proof} We 
induct on the height of $\sigma$. If the height is 1, the sequence is
constant. For the induction step, replace $\sigma$ by the iterate
$[f^M(\sigma)]$ from Theorem \ref{improved} so that it splits into
subpaths which are either single edges or exceptional paths. It suffices
to prove the statement for these subpaths. The statement is
clear for the exceptional subpaths, and it follows from the inductive
assumption for single edges.\end{proof}

More generally, we can consider polynomial sequences in any $\f$-tree.

\definition Let $T$ be an $\f$-tree. A sequence $\{ A_k\}_{k=k_0}^{\infty}$ of
embedded paths in $T$ is said to be {\it polynomial} if
it can be obtained from constant sequences of paths by finitely many
operations described below. \begin{enumerate} 
\item[0.] ({\it translation}): $A_k=\gamma_k(B_k)$ for a polynomial
sequence
$\{ B_k\}$ and a sequence $\{ \gamma_k\}$ of elements of
$\f$,
\item[1.] ({\it reindexing and truncation}): $A_k=B_{k+k'}$ for a
polynomial sequence $\{ B_k\}_{k=k_1}^{\infty}$ for
some $k_1\leq k_0+k'$, 
\item[2.] ({\it inversion}): $A_k$ is the inverse of $B_k$, and $\{
B_k\}_{k=k_0}^{\infty}$ is a polynomial sequence, \item[3.] ({\it
concatenation}): $A_k=B_k C_k$, where $\{ B_k\}$ and $\{ C_k\}$ are
polynomial sequences, (and in particular no cancellation occurs in
$B_kC_k$), and \item[4.] ({\it integration}): $A_k=B_{k_0} B_{k_0+1}\cdots
B_k$, where $\{ B_k\}$ is a polynomial sequence (and again there is no
cancellation). \endroster\enddefinition
 
The following lemma is by induction on the number of above operations
and its proof is left to the reader.
 
\begin{lem} \label{paths} Let $\{ A_k\} _{k=k_0}^{\infty}$ be a polynomial 
sequence of
paths in an $\f$-tree $T$. Then \roster
\item the function $k\mapsto \ell_T(A_k)$ is a polynomial function in $k$,
\item if $A_k=B C_k D$ for some paths $B$ and $D$, then the
sequence $\{ C_k\} _{k=k_0}^{\infty}$ is polynomial, 
\item $\{ A_k\}$ is either constant (up to the action of $\f$), or for
any $d>0$ there is $k_1\geq k_0$ so that for every $k\geq k_1$ $A_k=B
C_k D$ for paths $B$ and $D$ of length $\geq d$,
\item the initial endpoints of the $A_k$'s lie in a single $\f$-orbit, and
similarly the terminal endpoints lie in a single $\f$-orbit.\qed
\endroster\end{lem}
 
\begin{prop} \label{noperiodic}
If $\oa$ is a $UPG({\f})$ automorphism, then all $\oa$-periodic
conjugacy classes are fixed.\end{prop}

\begin{proof} Assume $x$ is an $\oa$-periodic conjugacy class. Represent
$x$ as a loop $\gamma$ in a $UR$ representative $f:\G\to
\G$. Consider the splitting of $\gamma$ given by the topmost edge of
$\G$ that intersects $\gamma$. Since $x$ is $\oa$-periodic, each of
the resulting paths is also $f$-periodic. Theorem \ref{eventually
polynomial} and Lemma \ref{paths} now imply that each path is
$f$-fixed, and thus $\gamma$ is $f$-fixed. \end{proof}

Another immediate consequence of Lemma \ref{paths} is the following.
 
\begin{prop} Let $\alpha :\tilde G\to T$ be an equivariant
map from the universal cover of $G$ to an $\f$-tree $T$ with a finite
$BCC$ (see Section \ref{bccsec}). Suppose $\{
A_k\}_{k=k_0}^{\infty}$ is a polynomial sequence in $\tilde G$, and
define $B_k=[\alpha (A_k)]$. Then
there is $k_1\geq k_0$ such that the sequence $\{
B_k\}_{k=k_1}^{\infty}$ is polynomial.\end{prop}
 
\begin{proof} By induction on the number of operations required to
construct $\{ A_k\}$. Focus on the last operation. Say $A_k=X_k
Y_k$. Inductively, we know that the sequences $\{ [\alpha (X_k)]\}$ and
$\{ [\alpha (Y_k)]\}$ are polynomial, after truncation. At most a bounded
amount can be canceled. Assuming they are
not constant, it follows from Lemma \ref{paths} (3) that eventually the
cancelled portions are independent of $k$, and the claim follows from
Lemma \ref{paths} (2). If one or both of the sequences are constant, the proof
is similar. 
 
The other nonobvious case (integration) is left to the reader.\end{proof}

\begin{thm} \label{limitexists} If $T\in\xx$ and $\oa\in UPG_{\f}$, then
the sequence $\{ T\oa^k\}$ converges to a tree $T\oa^\infty\in\xx$. 
Further, if $T\in \x$ then $T\oa^\infty\in\x$.\end{thm}
 
\begin{proof} Let $f:G\to G$ be a $UR$ for $\oa$.
If $w$ is any conjugacy class, the function $k\mapsto \ell_T(\oa^k(w))$
is eventually polynomial (the transition from paths to loops uses the
fact that $\tilde G\to T$ has a $BCC$ one more
time). The degree of the polynomial is uniformly bounded by the number
of strata in $G$. Let $d$ be the largest degree that occurs for this
$T$ and variable $w$. Then for any $w$ the sequence $k\mapsto
\ell_T(f^k(w))/k^d$ converges, and not all limits are $0$ thus the
sequence converges to an $\f$-tree. To see that the limiting tree is
simplicial if $T$ is simplicial, argue by induction on $d$ that if $\{
A_k\}$ is a polynomial sequence in $T$ (with endpoints in the vertex
set) of degree $m$, then the leading term of the polynomial $k\mapsto
\ell_T(A_k)$ is uniformly bounded away from 0. It follows that the
collection of nonzero numbers $\lim \ell_T(\oa^k(w))/k^d$ is bounded away from 0,
and so $T\oa^\infty$ is simplicial.

That this tree is very
small follows from the fact, proved in \cite{cl:verysmall},
that the collection of very small trees is closed under limits (in the
projectivized space). 
\end{proof}

\begin{defn} \label{growingdef}
Trees $T$ for which $d=0$ (i.e. the sequence $\{ \ell_T(\oa^k(w))\}$ is
eventually constant for every $w$) are called {\it
non-growers}. Others are {\it growers}. If $d=1$, we say that
$T$ {\it grows linearly}, etc.\end{defn} 

\begin{com} \label{nongrowers} 
There exist non-growers that are not fixed. An example is the tree $T$
with $T/<a,b>$ a circle with an arc attached at one endpoint, the
other endpoint labeled $<a>$, and all other labels trivial. The loop
corresponds to $b$, and $f (a)=a$, $f (b)=ab$. Such examples do not
exist in $SL_n(\mathbb Z)$ or the mapping class group of a
surface. Non-growers are also responsible for the existence of compact
sets $K\subset\xx$ in the complement of $Fix(f)$ with the property
that for no $k$ is $Kf ^k$ contained in a certain small neighborhood
of $Fix(f)$. A concrete example can be described as follows. Let
$F_4=<a,b_1,b_2,b_3>$, and let $f$ be given by $f (a)=a$ and $f
(b_i)=b_ia$. The compact set $K$ consists of the biinfinite sequence
$\dots,T_{-2},T_{-1},T_0,T_1,T_2,\dots$ together with the limiting
tree $T_{\infty}$. The quotient $T_n/F_4$ is the graph obtained from
the triod by attaching loops to the valence 1 points. The center point
is labeled $<a>$ and all other labels are 1. The three loops
correspond to $b_1a^n$, $b_2a^{-n}$, and $b_3$ respectively.  Passing
to the limit as $n\to \infty$ amounts to unfolding the first two loops
which in the limit correspond to $b_1$ and $b_2$ respectively. Now
notice that $T_nf^n$ converges to a non-fixed non-grower (which is a
tree just like $T_{\infty}$ except for a permutation of $\{
b_1,b_2,b_3\}$).
 
It is, however, true that if $K$ is a compact subset of the closure of
Outer Space consisting of growers, then the accumulation set of the
sequence $Kf^k$ is a subset of $Fix(f)$.
 
\end{com}

\subsection{Suffixes and eigenrays}

Recall that if $f:G\to G$ is a $UR$, an eigenray associated to an edge
$E_i$ is the infinite immersed path $E_iu_i[f(u_i)][f^2(u_i)]\cdots$
arising as the limit of iterates $[f^k(E_i)]$. The following proposition
is the analogue of the fact in linear algebra that if $A$ is a unipotent
matrix and $v$ a nonzero vector, then projectively the sequence $A^k(v)$
converges to an eigenspace of $A$.

\begin{prop} \label{more} Let $f:G\to G$ be a $UR$.
If $[f(u_i)]\not= u_i$, $R^*$ is an initial segment of $R_i$,
and
$\gamma$ is an immersed edge path in $G$ that contains
$E_i$ then there is an
$N$ such that, for all $k>N$, $[f^k(\gamma)]$ contains $R^*$ or its
inverse as a subpath.
\end{prop}

\begin{proof} 
We argue by induction on $height(\gamma)$. If $height(\gamma)=i$,
consider the splitting of $f^M(\gamma)$ into edges and exceptional paths.
There is
a 1-1 correspondence between occurrences of $E_i$ in $\gamma$ and in
$f^M(\gamma)$. Since
$[f(u_i)]\neq u_i$, $E_i$ does not occur in an exceptional path, and hence
one of the subpaths in the splitting is $E_i$ or $E_i^{-1}$.
Eventually, the iterates contain $R^*$ or its inverse.

Now assume $height(\gamma)=j>i$. Again consider the splitting of
$f^M(\gamma)$ into edges and exceptional paths.  First note that an
exceptional path $E_s\tau^k E_t^{-1}$ in this decomposition cannot
cross $E_i$ ($E_s$ and $E_t$ have fixed suffixes and so are distinct
from $E_i$, and $\tau$ cannot cross $E_i$ since $height(\tau)<j$ and
so otherwise by induction the iterates of $\tau$ (which equal $\tau$)
would have to contain arbitrarily long segments of $R_i$). If the edge
$E_i$ or its inverse occur in the splitting, we are done. Also, if
there is an edge $E_l$ in the splitting whose eigenray $R_l$ crosses
$E_i$, then large iterates of $\gamma$ contain large segments of
$R_l$, which in turn contain large iterates of $u_l$, and these
eventually contain $R^*$ by induction. It remains to exclude the
possibility that $E_i$ is not crossed by any of the eigenrays $R_l$ of
the edges $E_l$ in the splitting. The set of edges crossed by these
eigenrays union all edges with fixed suffixes is an $f$-invariant
subgraph (by induction) that contains $[f^m(\gamma)]$ for large $m$,
and does not contain $E_i$.  The restriction of $f$ to this subgraph
is a homotopy equivalence, and therefore $\gamma$ is homotopic into
it, contradicting the hypothesis.
\end{proof}

We next analyze the edge stabilizers of the tree obtained in the limit
under iteration by a $UPG$ automorphism, starting with certain
trees with trivial edge stabilizers that are closely related to a $UR$.
We discover that the edge stabilizers of the limiting tree are conjugates
of certain suffixes of the $UR$.

\begin{prop} \label{only u's pull}
Let $f:G\to G$ be a $UR$ of
$\oa\in Out(\f)$ and let $G_r$ be a subgraph in the associated filtration of
$G$. Assume that for every edge $E$ of $G$ we have $f(E)=Eu$ with $u$
either freely homotopic into $G_r$ or $[f(u)]=u$. Also assume that for at
least one such $u$ the first alternative fails.

Let $S$ be the tree obtained from the universal cover $\tilde G$ by
collapsing all edges that project into $G_r$. Then the stabilizer
of any edge in $T=Sf^{\infty}$ is infinite cyclic, and it contains
a conjugate of a nontrivial suffix that is not
freely homotopic into $G_r$.
\end{prop}

Example \ref{ex1} illustrates this phenomenon with $G_r=\emptyset$.

\begin{proof} Notice that $S$ grows linearly under $\oa$. Every path $\omega$
in $G$ determines a sequence of paths in $S$ by lifting the iterates
$[f^i(\omega)]$ to $\tilde G$ with a common initial point and
projecting to $S$. This sequence determines a path $f^\infty(\omega)$
in the limiting tree $T$ (thought of as the Gromov-Hausdorff limit of
the trees $S\oa^i$ scaled linearly, see \cite{fp:trees}). By
construction, paths of the form $f^\infty(\omega)$ cover $T$. Recall
that by Theorem \ref{improved} if $\omega$ is any path or a loop in
$G$, a sufficiently high iterate $[f^M(\omega)]$ has a splitting into
edges and exceptional subpaths. This splitting induces a subdivision
of $f^\infty(\omega)$ into subarcs, and shows that $T$ is covered
by paths of the form $f^\infty(\omega)$ where $\omega$ is an edge or
an exceptional path. If $\gamma$ fixes an arc in $T$, it must fix a subarc
of some such $f^\infty(\omega)$.

We next analyze the stabilizers of the nondegenerate arcs of the form
$f^\infty(\omega)$, and show that they each contain a
nontrivial suffix not homotopic into $G_r$.
Since $T$ is very small, any subarc of $f^\infty(\omega)$ has the same
stabilizer and the proposition follows.

Consider first an edge $E_i$ and the associated sequence
$E_iu_i[f(u_i)]\cdots [f^k(u_i)]$.  If $u_i$ is homotopic into $G_r$,
then the path $f^\infty(E_i)$ is degenerate, and if it is fixed, then
the path is fixed by $u_i$, viewed as an isometry of the limiting tree
(a typical element of the sequence is $E_i$ followed by a long string
of $u_i$'s, all contained in the axis of the isometry induced by
$u_i$, where we take the endpoint of $E_i$ as the basepoint). An
exceptional path $E_i\tau^k E_j^{-1}$ similarly determines a
degenerate path (if it is Nielsen) or a path fixed by $u_i$ (if it is
not). 
\end{proof}

\subsection {Primitive subgroups}
By looking at homology, it is clear that if free factors in a free
factor system are permuted under a $UPG$ automorphism, then they are
invariant, and the restriction is $UPG$.  We will show moreover that a
periodic free factor (or even a vertex stabilizer of a tree in $\xx$)
is invariant. Our argument uses only that vertex groups are primitive.

\begin{lem}\label{primitive} 
Let $H$ be a finitely generated primitive subgroup of ${\f}$, i.e. if
$\gamma^n\in H$ for some $n>0$ then $\gamma\in H$. Then the normalizer
$N(H)$ of $H$ in ${\f}$ is $H$.\end{lem}

\begin{proof} Let $T$ be a minimal free simplicial $\f$-tree and let $T_H$
be a minimal $H$-invariant subtree of $T$. Let
$\gamma\in N(H)$. Then
$\gamma(T_H)=T_H$ and so the axis of
$\gamma$ is in $T_H$ and projects to a loop in $T_H/H$. Thus, a power
of $\gamma$ is in $H$. Since $H$ is primitive, $\gamma$ is in
$H$.\end{proof}

We say that a subgroup $H$ of $\f$ is invariant under a subgroup ${\Cal
H}\subset Out(\f)$ if for every $\oa\in \Cal H$ and every lift $\hat \oa\in
Aut(\f)$ of
$\oa$, the subgroup $\hat \oa(H)$ is conjugate to $H$.

\begin{lem} \label{restriction} Let $\Cal H$ be a subgroup of $Out({\f})$
and let $H$ be a finitely generated primitive subgroup of ${\f}$ that
is $\Cal H$-invariant. Then the restriction map $\rho_H:{\Cal H}\to
Out(H)$ is well-defined. Further, if ${\Cal H}$ consists of $PG(\f)$
automorphisms, then ${\Cal H}|_H:=\rho_H({\Cal H})$ consists of
$PG(H)$ automorphisms.
\end{lem}

\begin{proof} The proof is an easy consequence of Lemma \ref{primitive}
and Lemma \ref{short}. \end{proof}

\begin{lem} \label{ifsomethenall}
Let $G'\to G$ be an immersion of finite graphs such that
$Im[\pi_1(G')\to\pi_1(G)]$ is a primitive finitely generated subgroup of
$\pi_1(G)$. Let
$\{ A_n\}$ be a polynomial sequence of paths in $G$. Assume that for infinitely
many values of
$n$ the path $A_n$ lifts to $G'$ starting at a given point $x\in G'$. Then the
same is true for all large $n$ and, furthermore, the lifts form (after
truncation) a polynomial sequence in $G'$ (so that in particular -- see Lemma
\ref{paths}(4) -- the terminal endpoint of these lifts is constant).
\end{lem}

The lemma fails if the primitivity assumption is dropped; e.g. take $G$ to be
the circle and $G'$ the double cover.

\begin{proof}
We proceed by induction on the number of basic operations in the construction
of $\{ A_n\}$. 

Suppose first that the last step is inversion. For infinitely many $n$ the
other endpoint of the lift of $A_n$ starting at $x$ is a point $y\in G'$
(there are finitely many preimages of the common terminal endpoint of the
$A_n$'s in $G$). Applying the statement of the lemma to $\{ A_n^{-1}\}$ we
learn that for all large $n$ there is a lift $\tilde A_n$ of $A_n$ that
terminates at $y$, for infinitely many $n$ it starts at $x$, and $\{\tilde
A_n$\} forms a polynomial sequence. Therefore, for all large $n$, $\tilde A_n$
starts at $x$.

Suppose next that the last step is concatenation:
$A_n=B_nC_n$. Then $B_n$ lifts to $G'$ starting at $x$ for infinitely many $n$
and thus for all large $n$, and the lifts $\tilde B_n$ form a polynomial
sequence. Let $y$ be the common terminal endpoint of the $\tilde B_n$.
Similarly, for all large $n$ the path $C_n$ lifts to a path $\tilde C_n$
starting at $y$, and these paths form a polynomial sequence. Thus $\tilde
A_n=\tilde B_n\tilde C_n$ is a polynomial sequence starting at $x$ and
projecting to $A_n$.

Finally, suppose that the last step is integration: $A_n=B_1B_2\cdots B_n$.
Since $A_n$ is a subpath of $A_{n+1}$ it follows from our assumptions that
each $A_n$ lifts to a path $\tilde A_n$ starting at $x_0=x$. Infinitely many
of these end at the same point $y_1$. Thus for infinitely many $n$ the path
$B_n$ lifts starting at $y_1$. It follows that eventually all these lifts end
at a point $y_2$. Again, for infinitely many $n$, $B_n$ lifts starting at
$y_2$ etc. Repeating this procedure we produce a sequence $y_1,y_2,\cdots$.
Suppose that $y_i=y_j$ for some $i<j$. For large $n$ there are lifts of $B_n$
that connect $y_i$ to $y_{i+1}$, $y_{i+1}$ to $y_{i+2}$,..., $y_{j-1}$ to
$y_j$. By the primitivity assumption we must have $y_i=y_{i+1}=\cdots=y_j$.
Therefore the sequence $y_1,y_2,\cdots$ is eventually constant, i.e. $y_n=y$
for all large $n$. Thus for large $n$ the path $B_n$ lifts to $\tilde B_n$
beginning and ending at $y$. The claim now follows.
\end{proof}

\begin{prop} \label{upgperiodic}
Suppose that $\oa$ is a $UPG({\f})$ automorphism and that $H\subseteq
\f$ is a primitive finitely generated subgroup. If $\oa^k(H)$ is
conjugate to $H$ for some $k>0$, then $\oa(H)$ is conjugate to
$H$. Furthermore, if $\hat \oa\in Aut(\f)$ is a lift of $\oa$ with $\hat
\oa^k(H)=H$ then $\hat \oa(H)=H$.
\end{prop}

The statement is false without the primitivity assumption as the
following example shows: $F_2=<a,b>$, $\hat\oa(a)=a$, $\hat\oa(b)=ab$,
$H=<a^2,b>$, $k=2$.

\begin{proof}
We may assume that $rank(H)>1$.
Let $f:G\to G$ be a $UR$ of $\oa$. By $p:\tilde
G\to G$ denote the covering space of $G$ corresponding to $H$. There is a lift
$F:\tilde G\to\tilde G$ of $f^k$. There is a
fixed point of $F$, perhaps after
replacing $F$ by a power. (Indeed, by linear algebra, some power $F^m$ of $F$
will have negative Lefschetz number. Any fixed point of negative index of 
$F^m$ composed with the retraction to the core is fixed under $F^m$.)
Let $v$ be a point fixed by $F$.

We now use $v$ and $p(v)$ as base points. Let $\alpha$ be a loop in $\tilde G$
based at $v$. The sequence $[F^i(\alpha)]$ of based loops forms a
sequence of lifts of a subsequence of the sequence $[f^j(p(\alpha))]$ of based
loops. The latter is eventually a polynomial sequence
(Theorem \ref{eventually polynomial}) and
hence by Lemma
\ref{ifsomethenall} for all large
$j$ the based loop
$[f^j(p(\alpha))]$ lifts to a based loop in $\tilde G$. Applying this to loops
$\alpha$ generating $\pi_1(\tilde G,v)$ we conclude that $f^j$ lifts to
$\tilde G$ for all large $j$. Thus $\oa^j(H)$ is conjugate to $H$ for large $j$
and the claim follows.

For the ``furthermore'' part of the proposition let $v$ be the base point
and choose $F$ so that $v$ is fixed.
\end{proof}

\begin{prop} \label{upgrestriction}
Suppose that $\oa$ is a $UPG({\f})$ automorphism and that $H\subseteq
\f$ is a primitive finitely generated subgroup. Then the restriction
(see Lemma \ref{restriction}) $\oa |_H$ is $UPG(H)$.
\end{prop}

\begin{proof}
Let $f:G\to G$ be a $UR$ of $\oa$. By $p:\tilde
G\to G$ denote the covering space of $G$ corresponding to $H$ and let $\tilde
f:\tilde G\to \tilde G$ be a lift of $f$. By $C$ denote the core of $\tilde
G$. Let
$\rho:\tilde G\to C$ be the nearest point retraction. If
$C$ does not contain any lifts of the topmost edge $E\subset G$, then we may
argue by induction on the number of strata. Therefore we assume that $C$
contains lifts of $E$. From $C$ form a finite graph $\G$ by collapsing all
complementary components of $C\setminus\cup\{ interiors\ of\ lifts\ of\ E\}$.
The map $\rho\tilde f$ induces a simplicial homeomorphism
$\phi:\G\to\G$. The main step of the proof is to argue that $\phi=id$.

Assuming $\phi\neq id$, we replace $f$ and $\phi$ by a power if
necessary so that there is a primitive loop $\gamma=E_1E_2\cdots E_m$
nontrivially rotated by $\phi$. Here each $E_i$ is a lift of $E$ or of
$E^{-1}$ and $\phi(E_i)=E_{i+r}$ for some $0<r<m$ (indices are $\mod
m$).

Choose a path in $C$ of the form $E_1\alpha E_2$ where $\alpha$ does
not cross any lifts of $E$ and denote by $\tau$ the subpath of
$E_1\alpha E_2$ obtained by splitting at $E_1$ and $E_2$. Now consider
the sequence $\{ [\tilde f^k(\tau)]\}_{k=1}^{\infty}$. This sequence
projects to an eventually polynomial sequence. Further, for infinitely
many values of $k$ (those in the same congruence class modulo the
order of $\phi$) these paths have common initial and common terminal
endpoints. It follows from Lemma \ref{ifsomethenall} that for large
$k$ and any $i$ there is a path that joins $E_{1+ir}$ and $E_{2+ir}$
and projects to the same path as $[\tilde f^k(\tau)]$. Repeat this
construction for every $\phi$-orbit of consecutive edges in $\G$ to
obtain a primitive loop in $C$ that projects to a proper power. This
contradicts the primitivity assumption and shows that $\phi=id$.

If $\alpha$ is any
loop in $C$ representing a cycle, then $(\tilde f_*-id)(\alpha)$ is a cycle
supported in the cores of the components of $C\setminus\cup\{ interiors\ of\
lifts\ of\ E\}$. Inductively, it follows that a high power of $\tilde f_*-id$
kills $\alpha$. Thus $\oa |_H$ is $UPG(H)$.
\end{proof}

\subsection{$UPG$ automorphisms and trees}

Recall that for $T\in\xx$ we denote by ${\cal V}(T)$ the set of point
stabilizers of $T$.  The set ${\cal E}(T)$ denotes the set of
stabilizers in $T$ of (nondegenerate) arcs. Note that $\cup {\cal
V}(T)$ is the set of elements of $F_n$ that are elliptic in $T$.

\begin{prop}\label{Eguts} Let $\oa$ be a $UPG({\f})$ automorphism, 
and let $T\in \xx$ such
that ${\Cal V}(T)$ and ${\Cal E}(T)$ are $\oa$-invariant. Then \roster
\item each element of $E\in {\Cal E}(T)$ is $\oa$-fixed (up to conjugacy),
\item each $V\in {\Cal V}(T)$ is $\oa$-invariant, and
\item the restriction of $\oa$ to each $V\in {\Cal V}(T)$ is $UPG(V)$.
\endroster
\end{prop}

\begin{proof}
The collection ${\Cal E}(T)$ consists of finitely many conjugacy
classes of cyclic subgroups of $\f$, by Theorem \ref{gl}. 
Therefore, generators of elements
of ${\Cal E}(T)$ are $\oa$-periodic, and hence $\oa$-fixed by Proposition
\ref{noperiodic}.

Similarly, each of finitely many representatives of conjugacy
classes in ${\Cal V}(T)$ is
$\oa$-periodic, hence $\oa$-invariant by Proposition \ref{upgperiodic}, and
the restriction of $\oa$ is $UPG$ by Proposition \ref{upgrestriction}.
\end{proof}

\begin{lem}\label{unifix} Let $g:G\to G$ be a simplicial homeomorphism 
of a connected finite graph. Suppose that $g$ fixes all valence one
vertices, and that either it induces identity map in $\hom(G)$ or that
it induces a unipotent map in $H_1(G;\mathbb Z)$. Then either $g=Id$
or $G$ is homeomorphic to $S^1$ and $g$ is rotation.
\end{lem}

\begin{proof} [Proof of Sublemma] First assume that $G$ has no valence
one vertices. If $G$ is a circle, the claim is clear. So assume $\chi
(G)<0$. By the Lefschetz fixed point theorem, $Fix(g)\neq
\emptyset$. Suppose $Fix(g)\neq G$.  Let $P$ be a shortest nontrivial
oriented edge path which intersects $Fix(g)$ only in its endpoints.

If any two $g$-iterates of $P$ either coincide or intersect only in
endpoints, then by considering the induced homomorphism on the homology of
the invariant subgraph of $G$ consisting of the union of all iterates of
$P$ we conclude that $g$ fixes $P$.

Suppose that there is an iterate $Q:=g^k(P)\neq P$ such that $P\cap Q$
contains a point that is not an endpoint of $P$. Then, unless $P\cap Q$ is
the common midpoint of $P$ and $Q$, $P$ is not the shortest nontrivial
oriented edge path which intersects $Fix(g)$ only in its endpoints. In
particular, $P\cap Q$ is fixed by $g^k$. So, replace
$g$ by
$g^k$, and $P$ by a proper subarc whose endpoints are fixed by $g^k$.
Repeating this will eventually construct a power of $g$ whose action on the
homology of a subgraph is not unipotent.

Now suppose that $G$ has a valence one vertex $v$. Any edge $E$ incident to $v$
must be $g$-fixed. So, remove $E$ and proceed by induction on the number of
edges.\end{proof}

Note that $\oa\in Out(F_n)$ fixes a very small tree $T$, i.e. if
$\ell_T(\oa(\gamma))=\ell_T(\gamma)$ for all $\gamma$, if and only if for any
lift $\hat\oa\in Aut(F_n)$ there is an $\hat\oa$-equivariant isometry
$f_{\hat\oa}:T\to T$.

\begin{prop} \label{directionfix} Assume $\n>1$. 
Suppose $\oa$ is a $UPG({\f})$ automorphism that fixes a
simplicial tree $T\in \x$. Let $\hat \oa\in Aut(F_n)$ be a lift
of $\oa$ and let $f_{\hat\oa}:T\to T$ be an $\hat\oa$-equivariant isometry.
Then $\oa$ fixes all orbits of vertices and
directions.\end{prop}

\begin{proof} The map $f_{\hat\oa}$ 
induces a periodic homeomorphism $\overline f_{\hat\oa}$ 
of the quotient graph.  It fixes all vertices whose labels are maximal
groups in $\cup\ve(T)$ by Proposition \ref{Eguts}. In particular, it
fixes all valence 1 vertices.  Since the induced action in homology of
the quotient graph is unipotent, by Lemma \ref{unifix}, $\overline
f_{\hat\oa}$ is identity or rotation of the circle. The latter is impossible
since then $\cup\ve(T)\neq \{ 1\}$.
\end{proof}

\begin{lem} \label{suffixgrows} 
Let $f:G\to G$ be a $UR$ for $\oa\in Out(F_n)$ and let $T\in \xx$.
Assume that whenever a suffix $u_i$ of $f$ is not fixed by $f$, then
there is a point in $T$ fixed by each $f^m(u_i)$, $m=0,1,2,\cdots$
(these are all loops based at the same point of $G$ and determine
elements of $\f$ up to simultaneous conjugacy). Then
\begin{enumerate}
\item $T$ is 
$\oa$-growing if and only if $\ell_T(u_i)>0$ for some suffix $u_i$,
and in that case the growth is linear.
\item Moreover, if $\ell_{T\oa^{\infty}}(\gamma)>0$ for a loop $\gamma$
in $G$ (see Theorem \ref{limitexists}), then there is a suffix $u_i$ as
in (1) such that for every
$N>0$ there exists $m_0>0$ with the property that all iterates
$[f^m(\gamma)]$, $m\geq m_0$, contain $[u_i^N]$ as a subpath.
\end{enumerate}
\end{lem}

\begin{proof} Let $\gamma$ be any loop in $G$. For large $m$, the loop
$[f^m(\gamma)]$ has a splitting $A_1(m)\cdot A_2(m)\cdot \cdots \cdot
A_k(m)$ into subpaths each of which is an edge or an exceptional path.
If there is an exceptional path $E_i\tau^k E_j^{-1}$ in this splitting
which is not Nielsen and with $\ell_T(\tau)>0$ then $u_i$ satisfies
(2). Similarly, if $u_i$ is the suffix associated to an edge in the
splitting with $\ell_T(u_i)>0$, then we have $[f(u_i)]=u_i$ by our
assumption, and again $u_i$ satisfies (2).  It remains to show that if
such $u_i$ does not exist, then $\ell_T(f^m(\gamma))$ remains bounded
as $m\to \infty$.

Let $\phi:\tilde G\to T$ be an equivariant map from the universal
cover of $G$ to $T$. For $l\geq m$ we have a splitting of
$[f^l(\gamma)]$ as $A_1(l)\cdot A_2(l)\cdot \cdots \cdot A_k(l)$
obtained by iterating the splitting above. Consider the lifts of these
paths to
$\tilde G$ starting at a fixed vertex $v$. Now argue inductively on $i$
that
$\phi$ sends the endpoint of the lift of $A_1(l)\cdot A_2(l)\cdot \cdots
\cdot A_i(l)$ to a point within bounded distance of $\phi(v)$. The inductive
step is clear when $A_i(l)$ is a Nielsen path. Now assume that $A_i(l)$
is an initial piece of an eigenray such that the associated suffix $u_i$
and all iterates of $u_i$ fix a point of $T$. First, the image of a lift
of the associated edge $E_i$ has finite length, so the image of its
endpoint is within bounded distance from $v$. The image of the
endpoint of $A_i(l)$ is obtained from this point by applying the group
element corresponding to a subpath of the eigenray $R_i$ that is a
concatenation of blocks.  This element is elliptic by assumption and
fixes a point of $T$ independent of $l$, so the claim follows.
Similar argument holds when $A_i(l)$ is exceptional or associated to
the inverse of an eigenray.
\end{proof}

\begin{lem} \label{eventually+} Let $\oa$ be a $UPG({\f})$ automorphism
and $S\in\xx$.  Suppose
$\ell_{S\oa^\infty }(\gamma)>0$. Then there is a $K$ such that
$\ell_{S}({\oa}^{K'}(\gamma))>0$ for all $K'>K$. \end{lem}
 
\begin{proof} This is immediate by the definition of limits.\end{proof}
 
If $f:G\to G$ is a homotopy equivalence that fixes all vertices of $G$ and
if $u$ is a path in $G$ with endpoints in the vertex set, then there
is a unique immersed path $f^{-1}(u)$ such that $[f(f^{-1}(u))]=[u]$.
 
\begin{prop}\label{no new 0's} 
Let $f:G\to G$ be a $UR$ for $\oa\in Out(F_n)$, and $T$ a tree in
$\xx$.  Assume that for each suffix $u_i$, a point in $T$ is fixed by
$u_i$, all its iterates $f^m([u_i])$ (so that $T$ is $\oa$-nongrowing),
and all negative iterates $f^{-m}([u_i^{-1}])$ of $[u_i^{-1}]$.  If
$\gamma$ is elliptic in $T\oa^{\infty}$, then $\gamma$ is elliptic in
$T$.
\end{prop}

\begin{proof} Represent $\gamma$ as a loop in $G$. Consider the splitting
of $[f^M(\gamma)]$ into edges and exceptional paths as
in Theorem \ref{improved}. We are assuming that when $m$ is
sufficiently large, $[f^m(\gamma)]$ lifts to a loop in the covering
space $G_V$ of $G$ corresponding to a vertex group $V$ of $T$.

It follows from Lemma \ref{ifsomethenall} that the splitting subpaths
of these lifts have endpoints independent of $m$ (for large $m$) and
also that the subpaths corresponding to blocks $[f^k(u_i)]$ are loops
for large $k$. We now claim that this is true for all $k$ and that paths
$[f^{-k}(u_i^{-1})]$ also lift to loops based at the same points.
Indeed, if $u_i$ is fixed by $f$, then there is nothing to prove, and if
it is not, then the group generated by $u_i$, its $f$-iterates and
$f^{-1}$-iterates of
$u_i^{-1}$ is nonabelian and fixes a unique point in $T$, and thus
is contained in a unique conjugate of $V$ which must be the one represented
by taking as basepoint the endpoints of the lift of $[f^k(u_i)]$ for
large $k$. Similarly, the subpaths of the lifts of $[f^m(\gamma)]$
corresponding to the exceptional paths $E_i\tau^k E_j^{-1}$ have the
property that $\tau$ forms a loop in $G_V$.

We now iterate $[f^m(\gamma)]$ backwards $m$ times and conclude that $\gamma$
lifts to a loop in $G_V$ since the lift is obtained from the lift
of $[f^m(\gamma)]$ by inserting loops of the form $f^{-k}(u_i^{-1})$.
\end{proof}

\section{A Kolchin Theorem for $UPG$ automorphisms}\label{kolchin}

The rest of the paper is devoted to the proof of our main theorem.

\begin{thm} \label{kolchin for trees} 
For every finitely generated $UPG(\f)$ group $\H$ there is a 
tree in $\x$ with all edge stabilizers trivial that is fixed by all elements
of $\H$ \end{thm}

\subsection{Bouncing sequences}

We start by setting up our iteration scheme, as outlined in the introduction.

\begin{defn} Let $\H$ be a $UPG$ group with a fixed finite generating set
$$\H=\langle \oa_1,\oa_2,\cdots,\oa_k\rangle$$ and let $T_0$ be any
simplicial tree in $\x$. The {\it bouncing sequence} associated with the above
data is the sequence of simplicial trees
$$T_0,T_1,T_2,\cdots$$ in $\x$ defined by 
$$T_{i}=T_{i-1}\oa_{i}^{\infty}$$ where subscripts of the $\oa_i$'s are
taken mod $k$ (see Theorem \ref{limitexists}).
\end{defn}

Notice that $T_i$ is $\oa_{i}$-fixed.
We will find a tree fixed by $\H$ by producing a bouncing sequence
that is eventually constant. In that case the stable value is a tree fixed
by all elements of $\H$.

\begin{ex} \label{ex1}
Let $F_2=\langle a,b\rangle$, $\H=\langle \oa\rangle$, with $\oa$ represented
by the automorphism $h:F_2\to F_2$ given by
$h(a)=a$, $h(b)=ba$,
and $T_0$ is a free simplicial $F_2$-tree. Then
$T_1=T_0\oa^{\infty}$ is the simplicial tree whose quotient graph
has one vertex labeled $\langle a,a^b\rangle$ and one edge (loop)
labeled $\langle a\rangle$. The loop is marked by $b$. This tree $T_1$ is
fixed by $\oa$ so the bouncing sequence is eventually
constant. However, $T_1$ has nontrivial edge stabilizers, and in this
case the iteration scheme fails to discover a tree as in the conclusion
of Theorem \ref{kolchin for trees}.\end{ex}

\begin{ex} \label{ex2} 
Let $F_3=\langle a,b,c\rangle$, $\H=\langle \oa_1,\oa_2\rangle$,
where $\oa_i$ is represented by $h_i$ given by 
$h_1(a)=a$, $h_1(b)=ba$, $h_1(c)=c$, $h_2(a)=a$, $h_2(b)=b$,
$h_2(c)=b^{-1}abc$.  Notice that the basis $\langle a,b,bc\rangle$ is
better adapted to $h_2$ since $h_2(bc)=abc$. Let $T_0$ be a simplicial
tree with trivial edge stabilizers whose quotient graph is the rose
with two petals marked $b$ and $bc$ respectively, and the single vertex
labeled $\langle a\rangle$.  Then $T_1=T_0\oa_1^{\infty}$ has
quotient graph a rose with petals marked $b$ and $c$, and the
vertex labeled $\langle a\rangle$. The tree $T_2=T_1\oa_2^{\infty}$ is
a tree combinatorially isomorphic to $T_0$, i.e.  $T_0$ and $T_2$
belong to the same simplex of $\x$. Notice, however, that $T_0$ and
$T_2$ are not homothetic: the ratio $length(b)/length(bc)$ is smaller
in $T_2$ than in $T_0$. The bouncing sequence indeed {\it bounces}
between two simplices in $\x$, so it does not stabilize. All trees in
the sequence are nongrowers under all elements of $\H$. The ratios
$length_{T_i}(b)/length_{T_i}(bc)$ converge to 0.\end{ex}

The above examples indicate the difficulties of trying to find 
a tree as in the conclusion of Theorem \ref{kolchin for trees} using
bouncing sequences. We will show, however, that the strategy is successful
provided we choose $T_0$ carefully.

\begin{thm} \label{bouncing ends well}
Let $\H=\langle \oa_1,\oa_2,\cdots,\oa_k\rangle$ be a $UPG$ group. By
$\Cal F$ denote a maximal $\H$-invariant proper free factor system.
Let $T_0$ be a simplicial tree with $\ve(T_0)=\Cal F$ and trivial edge
stabilizers. Then the bouncing sequence that starts with $T_0$ is
eventually constant, and the stable value is a simplicial tree with
trivial edge stabilizers.\end{thm}

In the beginning it is not clear that $\Cal F$ is nontrivial (although
this is a consequence of Theorem \ref{kolchin for trees}). The existence of
$\Cal F$ is guaranteed by Proposition \ref{ffs chain bound}.

In Example \ref{ex1} we started with $\Cal F$ trivial. The bouncing
sequence is eventually constant, but the edge stabilizers aren't
trivial. In this case we discover a larger invariant proper free
factor system, namely $\langle
a\rangle$ and its conjugates, by looking at the edge stabilizer.

In Example \ref{ex2} we started with $\Cal F$ consisting of $\langle
a\rangle$ and its conjugates. The sequence did not even
stabilize. However, we find a loop, namely $b$, that gets shorter and
shorter in the bouncing sequence (compared to other elements). This
tells us how to enlarge $\Cal F$ to a larger invariant free factor
system, namely $\langle a,b\rangle$ and its conjugates.

The proof of Theorem \ref{bouncing ends well} occupies the rest of
this section.

\subsection{Bouncing sequences grow at most linearly}

We now show that each tree in the bouncing sequence is either a nongrower,
or it grows linearly (assuming the choice of $T_0$ was made as in the
statement of Theorem \ref{bouncing ends well}).

\begin{prop} \label{ugly}
Let $\H$ be a $UPG$ group and let
${\Cal F}$ be an
$\H$-invariant proper free factor system of maximal complexity.
For
$\oa\in \H$, let
$f:G\to G$ be a $UR$ such that some subgraph $G_r$ in the
filtration of $G$ represents $\Cal F$. Let $E$ be an edge of $G$,
$u$ the corresponding suffix, and $R$ the corresponding eigenray, i.e.
$R=E\cdot u\cdot [f(u)]\cdot [f^2(u)]\cdots$.

Then at least one of the following holds.
\roster
\item The eigenray $R$ is eventually contained in $G_r$, or
\item $[f(u)]=u$.
\endroster
\end{prop}

Notice that $\Cal F$ is contained in
the vertex set of all trees in the bouncing sequence. Applying Lemma
\ref{suffixgrows} to a $UR$ $f_i:G_i\to G_i$ for $\oa_i$ 
we see that if all suffixes of $f_i$ are as
in (1), then $T_{i-1}$ is an $\oa_i$-nongrower, and otherwise it grows at
most linearly.

\begin{proof} 
Suppose the proposition fails for an edge $E$. We may assume that
$E$ is not crossed by any suffix of $f$, for if 
$f(E') = E'u'$ and $u'$ crosses $E$ then we may
replace $E$ by $E'$. Indeed, since $u\not= [f(u)]$, the eigenray $R' =
E'u'f(u')\cdots$ contains arbitrarily long subpaths of the eigenray
$R=Euf(u)\cdots$ by Proposition
\ref{more}. Thus
$R'$ crosses edges not in
$G_r$ infinitely often and it does not have periodic tail.

The edge $E$ determines a splitting of ${\f}$ as either a free product
or an HNN extension.  Let ${\Cal F}_E$ denote the resulting free
factor system. Note that $E$ is not an edge of $\G_r$ (otherwise $u$
would be in $\G_r$), and therefore ${\Cal F}\le {\Cal F}_E$. Also,
${\Cal F}\neq {\Cal F}_E$ since $R$ is contained in ${\Cal F}_E$ and it
is not eventually contained in $\Cal F$. 

The rest of the proof breaks into two cases. By $\hat \H$ denote the
preimage of $\H$ in $Aut({\f})$. Let $e$ be the point in $\partial \f$
determined by a lift of $R$ to the universal cover and let $\hat \H\{
R\}$ denote the set $\{ \{ \hat \oa e\} | \hat \oa\in \hat \H\}$ (this set
depends on $R$, but not on $e$).

{\bf Case 1:} $\hat \H\{ R\} \le \overline{{\Cal F}_E}$ (equivalently,
for all $\oa\in \H$ the ray $[\oa(R)]$ crosses $E$ only finitely many
times). In this case the smallest free factor system containing
$\Cal F$ and whose closure contains $\hat
\H\{ R\}$ (see Notation \ref{closure})) is proper (since it is contained in
${\Cal F}_E$),
$\H$-invariant (since both $\Cal F$ and $\hat \H\{ R\}$ are), and it
strictly contains $\Cal F$ (since $e$ is not in $\overline{\Cal F}$). This
contradicts the choice of $\Cal F$.

{\bf Case 2:} $\hat \H\{ R\} \not\le \overline{{\Cal F}_E}$. We will
show that in this case $\H$ contains an element of exponential
growth. There is $\oa\in \H$ such that, when represented as a homotopy
equivalence $g:G\to G$, $[g(R)]$ contains infinitely many $E$'s.  The
idea is that the image of a path containing $E$'s under a high power
of $f$ contains long initial subpaths of $R$ and the image under $g$
of a path with long initial subpaths of $R$ contains lots of
$E$'s. This feedback gives rise to exponential growth. We now make
this more precise. Let $R^*$ denote an initial subpath of $R$ chosen
long enough so that $[g(R^*)]$ contains 6 $E^{\pm 1}$'s with
occurrences of distance at least the $BCC$ (see Section \ref{bccsec})
for $g$ away from its endpoints.  Let $M$ be the length of $[g(R^*)]$.
Choose $N$ so that for all immersed paths $Ew$ and $EwE^{-1}$ where
$w$ is a path in $\G_r$ of length no more than $M$ we have that each
of $[f^N(Ew)]$ and $[f^N(EwE^{-1})]$ starts with $ER^*$. We claim that
the element of $\H$ represented by $gf^N$ has exponential growth.

Indeed, since $\f$ and the universal cover of $\G$ are quasiisometric,
it is enough to find a loop $\sigma$ in $\G$ such that the length of
$[(gf^N)^i(g(\sigma))]$ grows exponentially in $i$. We show that $\sigma$ can
be taken to be any
immersed based loop containing $ER^*$. In this case, $[g(\sigma)]$ contains
$[g(R^*)]$ except that perhaps subpaths containing endpoints of length
less than the $BCC$ for $g$ may have been
lost. In particular, $[g(\sigma)]$ contains 6 $E^{\pm 1}$'s separated by a
distance of no more than $M$. So it contains at least two disjoint
immersed subpaths of the form $(EwE^{\pm 1})^{\pm 1}$ where $w$ is a
path in $\G_r$ of length no more than $M$. Since $E$ is topmost, by
Proposition \ref{more}, $[f^Ng(\sigma)]$ contains two disjoint subpaths
of the form $(ER^*)^{\pm 1}$. So, $[gf^Ng(\sigma)]$ contains 2 disjoint
copies of $g(R^*)$ except for a loss of paths of length less than the
$BCC$ for $g$ and so contains at least 2
disjoint subpaths each with 6 $E^{\pm 1}$'s that are separated by a
distance of no more than $M$.  This pattern continues and the number
of such subpaths containing 6 $E^{\pm 1}$'s at least doubles with
application of $gf^N$.
\end{proof}

\subsection{Bouncing sequences stop growing}

Let $\oa\in Out({\f})$.
Recall from Definition \ref{growingdef} that a tree
$T\in\x$ is
$\oa$-growing if there is $\gamma\in{\f}$ such that
$lim_{m\to\infty}\ell_T(\oa^m([\gamma]))=\infty$.

\begin{prop}\label{evnotgrow} Let 
$\H=<\oa_1,\cdots,\oa_k>$ be a $UPG$ group, and let $T_0,T_1,\cdots$ be
a bouncing sequence for $\H$ as in Theorem \ref{bouncing ends
well}. Then all but finitely many elements of the sequence are
$\oa_i$-nongrowing for $i=1,2,\cdots,k$.
\end{prop}

\begin{proof} For notational simplicity, 
we assume that $\H=<\oa_1,\oa_2>$ and show that in the sequence $$T_0,
S_0:=T_0\oa_1^\infty, T_1 := S_0\oa_2^\infty, S_1 := T_1\oa_1^\infty,
T_2 := S_1\oa_2^\infty, \cdots$$ only finitely many elements are
$\oa_1$-growing. We will identify homotopy classes of elements of $\f$
with immersed loops in marked graphs.  Choose a $UR$ $f:G\to G$ for
$\oa_1$ so that $\cal F$ is represented by an invariant subgraph
$G_r$. Let $\cal U$ be the (finite) set of suffixes of $f$ that are
fixed by $f$. Set $K=|\cal U|$.  In fact, we will show that at most
$K$ of the $T_i$'s can be $\oa_1$-growing. Indeed, suppose that
$T_{i_0}, T_{i_1}, \cdots, T_{i_K}$ are $\oa_1$-growing with
$i_0<i_1<\cdots<i_K$.  By Lemma \ref{suffixgrows}, there is a suffix
$u_K$ of $f$ such that $\ell_{T_{i_K}}(u_K) > 0$. Thus $u_K$ (and its
$f$-iterates) are not elliptic in $T_0$ and in particular the eigenray
$$\cdots [f^s(u_K)]\cdot [f^{s+1}(u_K)]\cdot [f^{s+2}(u_K)]\cdots$$ is
not eventually contained in $G_r$. Therefore, by Proposition
\ref{ugly}, $u_K$ is fixed by $f$. Applying Lemma
\ref{eventually+} $2(i_K-i_{K-1})-1$ times, we see that there is a
word $w_K$ in $\oa_1$ and $\oa_2$ such that
$\ell_{S_{i_{K-1}}}(w_K(u_K))>0$. Lemma
\ref{suffixgrows} then provides a suffix $u_{K-1}$, such that
$\ell_{T_{i_{K-1}}}(u_{K-1}) > 0$ and, for large $B$,
$[f^B(w_K(u_K))]$ has a long string of $u_{K-1}$'s.  Continuing in
this fashion, we establish

\begin{sublem}\label{words} There are words $w_i \in <\oa_1,\oa_2>$, 
$1\le i \le
K$ and $u_i \in \cal U$, $0\le i\le K$ such that, for large $B$,
$[f^B(w_i(u_i))]$ contains a long string of $u_{i-1}$'s.\end{sublem}

Two of the $u_i$'s are equal, say $u_0=u_K$. We next find an element
in $<\oa_1,\oa_2>$ of exponential growth, a contradiction that will establish
the proposition.

Let $C$ be as in Theorem \ref{improved}(5) for the $UR$ $f$ and
choose $B$ so that the immersed based loop $[f^{B}w_i(u_i)]$ contains
$u_{i-1}^{C+2+A}$ where $A$ is chosen so that the length of $u_i^A$ is
larger than twice the maximum of the $BCC$'s of
the $w_i$'s (realized as homotopy equivalences on $G$). Then
$\oa_1^{B}w_1\dots \oa_1^{B}w_K$ has exponential growth. Indeed, we
will show that if $\gamma$ is any immersed path in $\G$ containing $L$
disjoint occurrences of $u_i^{C+2+A}$ then $[f^Bw_i(\gamma)]$ contains
$2L$ disjoint occurrences of $u_{i-1}^{C+2+A}$. After all, when we
apply $w_i$ to $\gamma$, we obtain for each occurrence of
$u_i^{C+2+A}$ an occurrence of $[w_{i}(u_{i}^{C+2})]$, the loss due to
the cancellation constant for $w_{i}$. So, by Theorem
\ref{improved}(5), each such occurrence gives rise to a splitting and,
upon application of $f^B$, we see $[f^Bw_{i}(u_{i}^{2})]$ which in
turn contains two disjoint copies of $u_{i-1}^{C+2+A}$. This ends the
proof of Proposition \ref{evnotgrow}.\end{proof}

\subsection{Edge stabilizers are eventually trivial}

We need the following lemma. Recall that for us an arc in a tree is a
subset homeomorphic to $[0,1]$.

\begin{lem} \label{edges of nongrowers}
Suppose that $T$ is a tree in $\xx$, $\oa$ is a $UPG$ automorphism, and $T$ is
$\oa$-nongrowing. Then every arc stabilizer of $T'=T\oa^{\infty}$ also
stabilizes an arc of $T$ and it is $\oa$-invariant.\end{lem}

\begin{proof}
Let $E=<e>$ be a nontrivial arc stabilizer of $T'$. Find an arc
$[v,w]$ in $T'$ that has an arc in common with $Fix_T(E)$ and two
elliptics $x$ and $y$ such that $Fix_{T'}(x)\cap [v,w]=\{ v\}$ and
$Fix_{T'}(y)\cap [v,w]=\{ w\}$.  Then we have that $x,y,e$ are
elliptics in $T'$ and
$\ell_{T'}(xy)>\ell_{T'}(xe)+\ell_{T'}(ye)$. Since $T$ is
$\oa$-nongrowing, for large $m$ we have
$\ell_{T'}(xy)=\ell_T(\oa^m(xy))$, etc.  Therefore,
$\ell_{T}(\oa^m(xy))>\ell_{T}(\oa^m(xe))+\ell_{T}(\oa^m(ye))$, and
$\oa^m(x)$, $\oa^m(y)$, $\oa^m(e)$ are elliptics in $T$. Hence
$\oa^m(e)$ is an edge stabilizer of $T$ for all large $m$. Since there
are only finitely many conjugacy classes of edge stabilizers in $T$,
it follows that the sequence $\oa^m(e)$ takes only finitely many
values, and is therefore constant (up to conjugacy) by Proposition
\ref{noperiodic}, and the lemma follows.\end{proof}

\begin{prop} \label{no z's}
The bouncing sequence $T_0,T_1,\cdots$ for
$\H$ in Theorem \ref{bouncing ends well} eventually consists of
trees that are $\oa_i$-nongrowing for all $i$ and have trivial edge
stabilizers. Further, for large $j$, the vertex stabilizers of $T_j$
are $\H$-invariant and independent of $j$.
\end{prop}

\begin{proof}
Eventually, the sequence consists of nongrowers by Proposition \ref{evnotgrow}.
Thus, eventually, the collection $\cup\ve(T_i)$ of elliptics forms a
nonincreasing sequence, by Proposition \ref{no new 0's}. It follows from
Proposition \ref{propC} that
eventually the sequence $\cup\ve(T_i)$ stabilizes. By Lemma \ref{edges
of nongrowers} eventually the collection of edge stabilizers stabilizes
as well. Let $T=T_j$ for some large $j$. Then $\cup\ve(T)$ is
$\Cal H$-invariant and contains $\cup\Cal F$, and all edge stabilizers of
$T$ are $\H$-invariant.

It remains to show that all edge stabilizers of $T$ are
trivial. Suppose $E$ is a nontrivial edge stabilizer of $T$. Let $p$
be the smallest integer such that $E$ fixes an edge of $T_p$. By our
choice of $T_0$, $p>0$. Lemma \ref{edges of nongrowers} implies that
$T_{p-1}$ is an $\oa_p$-grower (subscripts of $\oa_i$'s are taken mod
$k$). We now apply Lemma \ref{only u's pull} to a $UR$ $f:G\to G$ for
$\oa_p$ and with $G_r$ corresponding to $\Cal F$. Since $\cup{\Cal
F}\subseteq \cup\ve(T_{p-1})$, there is an equivariant map $\phi:S\to
T_{p-1}$, where the tree $S$ is obtained from the universal cover of
$G$ by collapsing all edges that project into $G_r$ as in Proposition
\ref{only u's pull}. In particular, there is a suffix of $f$ that is
not elliptic in $S$, so the hypotheses of Proposition \ref{only u's pull}
are satisfied. Thus both $S$ and $T_{p-1}$ grow linearly under
$\oa_{p}$. The map $\phi$ has finite $BCC$ (by Proposition
\ref{splittech}) Therefore the $BCC$ of the induced equivariant map
between $S\oa_{p}^m$ and $T_{p-1}\oa_{p}^m$, after scaling by $1/m$,
converges to 0 as $m\to \infty$. In the limit we obtain an equivariant
map $S\oa_{p}^{\infty}\to T_{p-1}\oa_{p}^{\infty}=T_p$ with
$BCC=0$. We conclude that $T_p$ is obtained from $S\oa_{p}^{\infty}$
by collapsing some edges and changing the metric on others. In
particular, $E$ fixes an edge of $S\oa_{p}^{\infty}$. By Proposition
\ref{only u's pull}, $E$ contains a conjugate of a suffix of $f$ not
homotopic into $G_r$. Now note that the free factor system given by a
topmost edge of $G$ contains both $\Cal F$ and $E$. Therefore, the
smallest free factor system that contains both $\Cal F$ and $E$ is
proper, and it is also $\H$-invariant (since $\Cal F$ and $E$ are),
and it properly contains $\Cal F$ (since it contains $E$, while $\Cal
F$ doesn't). This contradicts the choice of $\Cal F$.
\end{proof}

\subsection{Finding Nielsen pairs}

\begin{defn} \label{nielsen pairs}
Let $T$ be a simplicial $\f$-tree with all edge stabilizers trivial, and
let
$\H$ be a $UPG$ group. Assume that all vertex stabilizers of $T$ are
$\oa$-invariant (up to conjugacy) for all $\oa\in \H$. We say that two
distinct nontrivial vertex stabilizers
$V$ and
$W$ of $T$ form a {\it Nielsen pair} for $\H$ if for all $\oa\in\H$ and all
lifts $\hat \oa$ of
$\oa$ to $Aut(\f)$ there exists $\gamma\in\f$ such that $\hat
\oa(V)=V^{\gamma}$ and $\hat \oa(W)=W^{\gamma}$. (It suffices to
check this for one lift.)\end{defn}

For example, if $T$ is fixed by $\H$ and $V$, $W$ are nontrivial stabilizers 
of neighboring vertices, then $V$ and $W$ form a Nielsen pair. 

The proof of the following facts is left to the reader.

\begin{lem} \label{basic facts}
Let $T$ and $\H$ be as in Definition
\ref{nielsen pairs}.
\begin{itemize}
\item If $T'$ is another simplicial $\f$-tree that has the same vertex
stabilizers as $T$, then two vertex stabilizers $V$ and $W$ form a
Nielsen pair in $T$ if and only if they form a Nielsen pair in $T'$.
\item If $\H=\langle \oa_1,\oa_2,\cdots,\oa_k\rangle$ and two vertex
stabilizers $V$ and $W$ of $T$ form a Nielsen pair for $\langle
\oa_i\rangle$ for all $i$, then they form a Nielsen pair for $\H$.
\end{itemize}\qed
\end{lem}

\begin{prop}\label{nielsen} Let $\H=<\oa_1,\cdots,\oa_k>$ be a $UPG({\f})$
group and let $T$ be a simplicial tree such that
\begin{itemize}
\item $T$ has trivial edge stabilizers,
\item $\ve(T)$ is $\H$-invariant, and
\item $T$ is $\oa_i$-nongrowing for all $i$.
\end{itemize}
Then $T$ contains a Nielsen pair for $\H$. \end{prop}

By $h_i:G_i\to G_i$ denote an $RTT$ representative of $\oa_i$ with an
invariant subgraph $G_i'$ corresponding to $\ve(T)$, and whenever $E$
is an edge outside $G_i'$, then $h_i(E)=uEv$ for closed paths $u$ and
$v$ in $G_i'$. Such a representative can be constructed from
$T\oa_i^{\infty}$ (which is a tree with the same set of elliptics as
$T$ by Proposition \ref{no new 0's}, but is $\oa_i$-fixed) by passing
to the quotient and blowing up vertices to $UR$'s of the restriction
maps. As usual, the indices of $h_i$'s and $\oa_i$'s are taken $mod\
k$. Using Lemma \ref{basic facts} we shall detect that two vertex
stabilizers $V$ and $W$ of $T$ form a Nielsen pair for $\H$ by
examining for every $i$ whether they form a Nielsen pair for $\langle
\oa_i\rangle$ in the tree $T_i$ obtained from the universal cover of
$G_i$ by collapsing all edges that project to $G_i'$.

Edge paths $P$ in $G_i$ are of the form $v_0P_1v_1P_2\dots P_pv_p$
where each $P_j$ is an edge not in $G_i'$ and each $v_j$ is a
path in $G_i'$.  We call the elements $v_j$ {\it vertex elements}
(referring to the vertices of $T$). Some of the $v_j$'s could be
trivial paths. When
$P$ is such a path, then the iterates $h_i^N(P)$ have a similar form 
$v_0^{(N)}P_1v_1^{(N)}P_2\dots P_pv_p^{(N)}$. For each $j$ the
sequence $v_j^{(N)}$ is eventually polynomial. We say that the
vertex element $v_j$ is {\it inactive} if $v_j^{(N)}$ is independent
of
$N$. Otherwise,
$v_j$ is {\it active}. Of course, $h_i$ and the edge path $P$ are
implicit in these definitions. Even trivial $v_j$'s could be
active. 

When $i\neq j$ there is a homotopy equivalence $\phi_{ij}:G_i\to G_j$
given by markings. We may assume that this map sends
vertices to vertices and restricts to a homotopy equivalence $G_i'\to
G_j'$. Let $C$ be a constant larger than the $BCC$
of any $\phi_{ij}$.
Let
$v$ be a vertex element in a path $P$ in $G_i$. We can transfer $P$ to
another
$G_j$ using
$\phi_{ij}$ and tightening. The path $\phi_{ij}(v)$ has length
bounded above and below by a linear function in the length of $v$,
and then at most $2C$ is added or subtracted due to the $BCC$. 
In particular, if the length of a vertex
element in $P$ is larger than some constant $C_0>2C$, then this vertex
element {\it induces} a well-defined vertex element in $G_j$. Short
vertex elements in $P$ can disappear and new short vertex elements
can appear in $[\phi_{ij}(P)]$.

Choose constants $C_1,C_2,\cdots,C_{7k}$ such that if a vertex element
$v$ has length $\leq C_i$ and is transferred to some other graph, then
the induced vertex element has length $\leq C_{i+1}$.
Also, fix
$\e\in (0,1/14k)$.

\begin{lem} \label{intervals}
For a sufficiently large integer $m>0$ the
following statements hold.

\begin{itemize}
\item Let $N_i = 2^{2^{(7k-i+1)m}}$,
and let $I_{i,l}$ be the interval $$\big[ (1-l\e)N_i,
(1+l\e)N_i^m\big]$$ for $i=1,2,\cdots,7k$, $l=1,2,\cdots,14k$. Then
$I_{i,1}\subset I_{i,2}\subset \cdots
\subset I_{i,14k}$ and the intervals $I_{i,14k}$ are pairwise disjoint
for $i=1,2,\cdots,7k$, and further, they are disjoint from $[0,C_{7k}]$.

\item If a vertex element $v$ in an edge path $P$ in
$G_i$ is active and has length $\leq (1+14k\e)N_{i+1}^m$ (which is
the right-hand endpoint of $I_{i+1,14k}$), then the
$h_i$-iterated vertex element $v^{(N_i)}$ has length in $I_{i,1}$.

\item If a vertex element $v$ in an edge path $P$ in
$G_i$ has length in $I_{p,l}$ ($l<14k$), then after transferring to
$G_j$
$v$ induces a vertex element whose length belongs to
$I_{p,l+1}$.

\item If a vertex element $v$ in an edge path $P$ in $G_i$
has length in $I_{j,l}$ and if $i>j$ and $l<14k$, then the iterated
vertex element $v^{(N_i)}$ in $h_i^{N_i}(P)$ has length in
$I_{j,l+1}$.
\end{itemize}
\end{lem}

We think of the first index in intervals $I_{i,l}$ as measuring the
order of magnitude of lengths of vertex elements. The second index is
present only for technical reasons: there is a slight loss when
transferring from one graph to another (bullet 3), and when applying
``lower magnitude maps'' (bullet 4).

\begin{proof}[Proof of Lemma \ref{intervals}]
To see that the right-hand endpoint of $I_{i+1,14k}$ is to the left of
the left-hand endpoint of $I_{i,14k}$ we have to show that
$$(1+14k\e)2^{2^{(7k-i)m}+m}<(1-14k\e)2^{2^{(7k-i+1)m}}$$
i.e. that
$$2^{[2^{(7k-i+1)m}-2^{(7k-i)m}-m]}>\frac{1+14k\e}{1-14k\e}$$
That the latter inequality holds for large $m$ follows from the
observation that the exponent of the left-hand side
$$2^{(7k-i)m}(2^m-1)-m$$
goes to infinity as $m\to \infty$.

It follows from Theorem \ref{improved}(4) that there are polynomials
$Q_i$ and
$R_i$ with nonnegative coefficients such that whenever $v$ is an active
vertex element in a path
$P$ in $G_i$, then the length of
$v^{(N)}$ is in the interval
$[N-R_i(|v|),(1+|v|)Q_i(N)]$. The proof now reduces to the fact that
exponential functions grow faster than polynomial functions.
For example, the second bullet amounts to the inequalities
$$N_i-R_i((1+14k\e)N_{i+1}^m)>(1-14k\e)N_i$$ and
$$(1+(1+14k\e))N_{i+1}^mQ_i(N_i)<(1+14k\e)N_i^m$$
If we assume without loss of generality that $R_i(x)=x^d$ then the
first inequality simplifies to 
$$\frac{N_i}{N_{i+1}^{m+d}}>\frac{(1+14k\e)^d}{14k\e}$$
Again, the left-hand side amounts to $2^{exp}$ with
$$exp=2^{(7k-i)m}(2^m-m-d)$$
and it goes to infinity as $m\to\infty$. The proof of the second
inequality and of the other claims in the lemma are similar. (For the
third bullet use the fact that there is a linear function $L$ such
that if $w$ is a vertex element of a path $P'$ induced by a vertex
element $v$ of a path $P$, then the length of $w$ is bounded by
$L(|v|)$.)
\end{proof}

We will argue that if there are no $\H$-Nielsen pairs in $T$,
then the element $\oa_{7k}^{N_{7k}}\cdots \oa_2^{N_2}\oa_1^{N_1}\in\H$ has
exponential growth.

Start with an immersed loop $P_1$ in $G_1$ that is not contained in
$G_1'$ and all of whose vertex elements have length $\leq C_1$. This
loop is the {\it first generation}.  Then apply $h_1^{N_1}$ to obtain
$h_1^{N_1}(P_1)$ and transfer this new loop via $\phi_{12}$ to
$G_2$. The resulting loop $P_2$ is the {\it second generation}.  Then
apply $h_2^{N_2}$ and transfer to $G_3$ to obtain the {\it third
generation} loop $P_3$ etc. The loop $P_{7k}$ whose generation is $7k$
lives in $G_{7k}$. Then repeat this process cyclically: apply
$h_{7k}^{N_{7k}}$ and transfer to $G_1$ to get a loop $P_{7k+1}$ of
$(7k+1)^{st}$ generation etc.

Suppose that $v$ is a vertex element of some $P_i$. If $v^{(N_i)}$ has length
$\geq C_0$, then $v^{(N_i)}$ induces a well-defined vertex element $v'$
in $P_{i+1}$. We say that $v$ {\it gives rise} to $v'$.

We will now label some of the vertex elements of the $P_i$'s with
positive integers.  Consider maximal (finite or infinite) chains
$u_1,u_2,\cdots$ of vertex elements such that $u_i$ gives rise to
$u_{i+1}$. In particular, there is an integer $s$ such that $u_i$ is a
vertex element of $P_{i+s}$ for $i\geq 1$. If the length of the chain is
$\geq 7k$, then label $u_i$ by the integer $i$. If the chain has $<7k$
vertex elements, we will leave all of them unlabeled. All labels
$>1$ in $P_i$ correspond to unique labels in $P_{i-1}$. A {\it birth} is the
introduction of label 1. A {\it death} is an occurrence of a labeled vertex
element that does not give rise to any vertex elements in the next
generation. Any labeled vertex element can be traced backwards to
its birth. Traced forward, any labeled vertex element either eventually dies,
or lives forever (and the corresponding label goes to infinity).

\begin{lem} If a vertex element $v$ in some $P_i$ is not labeled, then
$v$ is $h_i$-inactive and its length is $\leq C_{7k}$.\end{lem}

\begin{proof}
The first element $v_1$ of a maximal chain $v_1,v_2,\cdots,v_s$,
$s<7k$, must have length $\leq C_1$.  Indeed, assume not. Say $v_1$ is
a vertex element in $P_{i+1}$. By the choice of $P_1$ we must have
$i\geq 1$. Transferring to $G_i$ $v_1$ induces a vertex element $v'$
of length $> C_0$. Now $v'=w^{(N_i)}$ and $w$ gives rise to $v_1$, so
the chain wasn't maximal.

If all $v_i$'s are inactive, then the claim about the length follows
from the definition of constants $C_i$. If $v_i$ is the first active
element of the chain, then $v_{i+1}$ has length in $I_{i,2}$ by the
second bullet of Lemma \ref{intervals}. With each generation the second
index of the interval increases by two until $7k$ generations are
complete (by bullets 3 and 4) or its length increases in length to
some
$I_{j,2}$ with $j<i$ by Property 2 and its life continues at least
$7k$ more generations. This  contradicts
$s<7k$.\end{proof}

\begin{lem} \label{event}
If two vertex elements in $P_i$ are labeled
with no labeled vertex elements between them, then either at least one
dies in the next $<k$ generations, or a birth occurs between them in
the next
$<k$ generations.
\end{lem}

\begin{proof}
If not, then the path between two such vertex elements is a Nielsen
path (i.e. its lift to $T$ connects two vertices whose stabilizers form
a Nielsen pair).\end{proof}

\begin{lem} \label{large labels}
Consider the cyclically ordered set of
labels in each
$P_i$.

\begin{itemize}

\item If two labels are adjacent, at least one is $<3k$.
\item If two labels have one label between them, then at least
one is $<4k$.
\item If two labels have two labels between them, then at least
one is $<5k$.
\item If two labels have three labels between them, then at least
one is $<6k$.
\end{itemize}
\end{lem}

\begin{proof}
Let $a$ and $b$ be two adjacent labels in some $P_i$ with
$a,b\geq 3k$ and assume that $i$ is the smallest such $i$. Consider the
ancestors of the two labels. According to Lemma
\ref{event} a death must occur between the two in some $P_{i-s}$ with
$s<k$. Thus in $P_{i-s}$ we have labels $\cdots (a-s)\cdots x\cdots
(b-s)\cdots$ and $x\geq 7k$. The dots between $(a-s)$ and $(b-s)$ are
vertex elements that die before reaching $P_i$, and their labels are
therefore $\geq 6k$. By our choice of $i$ we conclude that $x$ is the
only label between $(a-s)$ and $(b-s)$. Now consider further ancestors
of
$(a-s)$, $x$, and $(b-s)$. Again by Lemma \ref{event} a death must
occur between vertex elements labeled $(a-s)$ and $x$ in some
$P_{i-s-t}$ with $t<k$. We thus have two adjacent labels $\geq 5k$ in
$P_{i-s-t}$, contradicting the choice of $i$.

Now suppose that in some $P_i$ we have labels $\cdots axb\cdots$ and
$a,b\geq 4k$. By the first bullet we must have $x<3k$. If a death
occurs between $a$ and $x$, or between $b$ and $x$, in the previous
$k$ generations, then we obtain a contradiction to the first bullet.
If not, then by Lemma \ref{event} we conclude that $x<k$ and then we
have adjacent labels $a-x-1$ and $b-x-1$ in $P_{i-x-1}$
contradicting the first bullet.

Proofs of the last two bullets are analogous.
\end{proof} 

\begin{proof}[Proof of Proposition \ref{nielsen}]
Suppose that there are no $\H$-Nielsen pairs in $T$. Let
$C_0,C_1,\cdots,C_{7k}$ and $\e$ be constants as explained above.
Let $m$ be an integer satisfying Lemma \ref{intervals}, and consider
the labeling of vertex elements in paths $P_i$ as above. The fact that
$\oa_{7k}^{N_{7k}}\cdots \oa_2^{N_2}\oa_1^{N_1}\in \H$ grows exponentially
now follows from the observation that the number of labels in
$P_{i+k}$ is at least equal to the number of labels in $P_i$
multiplied by $5/4$. Indeed, consider the labels in $P_i$ that will
die before reaching
$P_{i+k}$. All such labels have to be $\geq 6k$ (since a vertex
element cannot die before reaching the ripe old age of
$7k$). By Lemma \ref{large labels}, any two such labels have at least 3
labels
$a$,
$b$, and
$c$ between them. By Lemma \ref{event}, there will be at least one
birth between
$a$ and $b$ and at least one birth between $b$ and $c$ between
generations
$i+1$ and $i+k$. Thus the number of deaths is at most a quarter of the number
of labels in $P_i$, and the number of births is at least twice the number of
deaths. The above inequality follows.
\end{proof}

\subsection{Distances between the vertices}

Consider the bouncing sequence as in Theorem \ref{bouncing ends well}.
Eventually, for $j\geq j_0$, $T_j$ is $\oa_i$-nongrowing for
$i=1,2,\cdots,k$ and the vertex groups of $T_j$ are $\H$-invariant. In
particular, the collection of vertex stabilizers of $T_j$ does not
depend on $j$. For $j\geq j_0$ we define the metric on
$T_{j+1}=T_j\oa_{j+1}^{\infty}$ by
$\ell_{T_{j+1}}(\gamma)=\ell_{T_j}(\oa_{j+1}^N(\gamma))$ for large $N$
(that is, we are taking the limit in the {\it unprojectivized} space
of trees). By Proposition \ref{nielsen} there is an $\H$-Nielsen pair
in $T_j$ for $j\geq j_0$.

\begin{lem} \label{distance}
Let $V$ and $W$ be two vertex stabilizers of $T_{j_0}$ and let
$d_j$ denote the distance between the vertices in $T_j$ fixed by $V$ and
$W$. If $V$ and $W$ form a Nielsen pair for $\H$, then $d_{j_0}=
d_{j_0+1}= d_{j_0+2}=\cdots$.
\end{lem}

\begin{proof}
Choose nontrivial elements $v\in V$ and $w\in W$. The distance between
the vertices in $T_j$ fixed by $V$ and $W$ equals
$\frac{1}{2}\ell_{T_j}(vw)$ and the distance in $T_{j+1}$ is analogously 
$\frac{1}{2}\ell_{T_{j+1}}(vw)$. The latter number can be computed as
$\frac{1}{2}\ell_{T_j}(\hat \oa_{j+1}^N(v)\hat \oa_{j+1}^N(w))$ for large
$N$, where $\hat \oa_{j+1}$ denotes a lift of $\oa_{j+1}$ to $Aut(\f)$ (since
$T_j$ is $\oa_{j+1}$-nongrowing). This in turn equals the distance in
$T_j$ between the vertices fixed by $\hat \oa_{j+1}^N(V)$ and $\hat
\oa_{j+1}^N(W)$. But that equals the distance between the vertices fixed by
$V$ and $W$ since $V$ and $W$ form a Nielsen pair for $\langle
\oa_{j+1}\rangle $. \end{proof}

\begin{lem} \label{nielsen edge}
Let $D_j\subset \mathbb R$ denote the set of distances between two
distinct vertices in $T_j$ with nontrivial stabilizer, $j\geq j_0$. Then
\begin{enumerate}
\item[(1)] $D_j$ is discrete,
\item[(2)] $D_j\supseteq D_{j+1}$ for all $j\geq j_0$, 
\item[(3)] there are finitely many $\f$-equivalence classes of paths $P$
joining two vertices of $T_j$ with nontrivial stabilizer and with
$length(P)=\min D_j$,
\item[(4)] if $V$ and $W$ are two nontrivial vertex stabilizers of $T_j$
such that the distance between the corresponding vertices is $\min D_j$,
then $V$ and $W$ form a Nielsen pair for $\langle \oa_j\rangle$.
\item[(5)] $\min D_j\leq \min D_{j+1}$, and
\item[(6)] if $\min D_j= \min D_{j+1}$ then any two nontrivial vertex
stabilizers $V$ and $W$ in $T_{j+1}$ realizing the minimal distance
also realize minimal distance in $T_j$.
\end{enumerate}
\end{lem}

\begin{proof}
(1) Every element of $D_j$ is a real number that can be represented as a
linear combination of (finitely many) edge lengths of $T_j$ with
nonnegative integer coefficients. Hence $D_j$ is discrete.

(2) Every element of $D_{j+1}$ has the form $\frac{1}{2}\ell_{T_j}(\hat
\oa_{j+1}^N(v)\hat \oa_{j+1}^N(w))$ (see the proof of Lemma \ref{distance})
and hence occurs also as an element of $D_j$.

(3) Let $P$ be such a path. The quotient map $T_j\to T_j/{\f}$ is
either injective on $P$ or identifies only the endpoints of $P$, hence
there are only finitely many possible images of $P$ in the quotient graph.
If two such paths have the same image, then they are
$\f$-equivalent.

(4) Since $\oa_j$ fixes $T_j$, for any lift $\hat \oa_j\in Aut(\f)$ of
$\oa_j$ we can choose an $\hat \oa_j$-invariant isometry $\phi:T_j\to
T_j$. By Proposition \ref{directionfix} and Lemma \ref{unifix}
$\phi$ induces identity in the quotient graph. Therefore the immersed
path $P$ joining the two vertices is mapped by $\phi$ to a translate
of itself (we are using the fact that all interior vertices of $P$
have trivial stabilizer).

(5) is a consequence of (2).

(6) Choose a lift $\hat \oa_{j+1}\in Aut(\f)$ of $\oa_{j+1}$. The
distance between the vertices corresponding to $V$ and $W$ has the form
$\frac{1}{2}\ell_{T_j}(\hat
\oa_{j+1}^N(v)\hat \oa_{j+1}^N(w))$ for large $N$. It follows that for large
$N$ the immersed path $P_N$ joining vertices in $T_j$ corresponding to
$\hat \oa_{j+1}^N(V)$ and $\hat
\oa_{j+1}^N(W)$ has length $\min D_j$. By (4), $V$ and $W$ form a Nielsen
pair for $h_{j+1}$ and therefore the paths $P_N$ are translates of each
other and have length $\min D_j$.
\end{proof}

\subsection{Proof of Theorem \ref{bouncing ends well}}

We are now ready for the proof of Theorem \ref{bouncing ends well}. For
the reader's convenience we first restate it.

\begin{bouncing kolchin}
Let $\H=\langle \oa_1,\oa_2,\cdots,\oa_k\rangle$ be a group in UPG. 
By $\Cal F$
denote a maximal $\H$-invariant proper free factor system. 
Let $T_0$ be a simplicial tree with
$\ve(T_0)=\Cal F$. Then the bouncing sequence that starts with $T_0$
is eventually constant, and the stable value is a simplicial tree with
trivial edge stabilizers.\end{bouncing kolchin}

\begin{proof}
The sequence eventually consists of nongrowers by Proposition
\ref{evnotgrow}. Then, eventually, the vertex stabilizers are
independent of the tree in the sequence and all edge stabilizers are
trivial by Proposition \ref{no z's}. By Proposition \ref{nielsen}
these advanced trees contain Nielsen pairs for $\H$. By Lemma
\ref{distance} it follows that the numbers $\min D_j$ of Lemma
\ref{nielsen edge} are bounded above and hence stabilize. Say $\min
D_{j+1}=\min D_{j+2}=\cdots=\min D_{j+k}$. Let $V$ and $W$ be two
nontrivial vertex stabilizers in $T_{j+k}$ that realize $\min
D_{j+k}$. By Lemma \ref{nielsen edge} $V$ and $W$ form a Nielsen pair
for every $\langle \oa_i\rangle$, and hence for $\H$. Let $P$ be the
immersed path joining the corresponding vertices. If $P$ projects onto
the quotient graph, then this quotient graph has one edge and
$T_{j+k}$ is fixed by $\H$. If $P$ does not project onto the quotient
graph, we obtain a contradiction by collapsing $P$ and its translates
and thus constructing an $\H$-invariant proper free factor system
strictly larger than $\cal F$.
\end{proof}

\section{Proof of the main theorem} \label{main proof}

In this section we show that Theorem \ref{kolchin for trees} implies
Theorem \ref{graph kolchin}. 

We start with an immediate consequence of Theorem \ref{kolchin for trees}.

\begin{prop} \label{lift} Every finitely generated $UPG$ group $\H$
lifts to a group
$\hat
\H\subset Aut(\f)$.\end{prop}

\begin{pf}
Let
$T$ be a simplicial
$\f$-tree with trivial edge stabilizers fixed by all elements of $\H$.
By collapsing orbits of edges we may assume that $T$ has only one orbit
of edges (the collapsing is possible by Proposition
\ref{directionfix}). Fix an edge
$e\subset T$. Since
$\oa\in\H$ fixes
$T$, there is a lift $\hat\oa\in Aut(\f)$ of $\oa$ and a
$\hat\oa$-equivariant isomorphism $f:T\to T$. We may choose
$\hat\oa$ and $f$ so that $f(e)=e$, and this choice is unique. The
set $\{\hat \oa|\oa\in\H\}$ is a group and gives the desired lift to
$Aut(\f)$.
\end{pf}

Recall from the introduction that for a
filtered marked graph $G$ the set of upper triangular homotopy
equivalences of $G$ up to homotopy relative to the vertices
is denoted by $\Q$.

\begin{lem} \label{gp str}  $\Q$ is a group under the
operation induced by
composition.
\end{lem}

\begin{pf} Since the
composition of
upper triangular homotopy equivalences  is clearly upper
triangular, it
 suffices to show that if  $\q$ is  upper triangular,
then there  exists an upper triangular $\r$ such that
$\q\r(E_i)$ and $\r\q(E_i)$ are homotopic rel endpoints to
$E_i$ for  $1
\le i \le K$.  We define
$\r(E_i)$ inductively starting with
$\r(E_1) = E_1$.  Assume that $\r$ is defined on $\ttt_{i-1}$
and that
$\q\r(E_j)$ and $\r\q(E_j)$ are homotopic rel endpoints to
$E_j$ for each
$j < i$.
 If $\q(E_i) = v_iE_iu_i$,  define $\r(E_i) = v_i'E_iu_i'$ where
$u_i'$ equals
$r(u_i)$ with its orientation reversed and $v_i'$ equals $r(v_i)$
with its
orientation reversed.  Since $v_i$ is a path in $\ttt_{i-1}$ with
endpoints at vertices, $\q\r(v_i)$ is homotopic rel endpoints to
$v_i$. Thus
$\q(v_i')$ is homotopic rel endpoints to $v_i$ with its orientation
reversed
and $v_i\q(v_i')$ is homotopic rel endpoints to the trivial path.
A similar
argument shows that $u_i\q(u_i')$ is homotopic rel endpoints to
the trivial
path and hence that  $\q\r(E_i) =
\q(v_i')v_iE_iu_i\q(u_i')$ is homotopic rel endpoints to $E_i$. A
similar
argument showing that $\r\q(E_i)$ is homotopic rel endpoints to
$E_i$
completes the proof.\end{pf}

\vskip .2cm

{\it Proof that Theorem \ref{kolchin for trees} implies Theorem
\ref{graph kolchin}}. Let $T$ be an $\H$-fixed 
tree with trivial edge stabilizers. As in the proof of Proposition
\ref{lift} we may assume that all edges of $T$ are translates of an
edge $e$. There are two cases depending on whether or not the
endpoints $a$ and $b$ of $e$ are in the same $\f$-orbit. We will
first consider the case that they are in distinct orbits, i.e. $T/\f$ is an
arc. By $A$ and $B$ denote the stabilizers of $a$ and $b$ respectively.
By induction on the rank, there exist desired
representatives $G_a$ and $G_b$ for $\H|_A$ and $\H|_B$
respectively. We define $G$ to be the disjoint union of $G_a$ and
$G_b$ with an edge $E$ connecting a vertex of $G_a$ and a vertex of
$G_b$. We choose a filtration of $G$ so that $E$ is the highest edge,
and so that this filtration induces the given once on $G_a$ and
$G_b$. For $\oa\in \H$ let $\hat\oa\in Aut(\f)$ and $f:T\to T$ be as
in the proof of Proposition
\ref{lift}. Let
$T_0$ be a free simplicial $\f$-tree and $f_0:T_0\to T_0$ a
$\hat\oa$-equivariant map. The triple $(T\times_{\f} T_0, (orbit\ of\
a)\times_{\f} T_0, (orbit\ of\ b)\times_{\f} T_0)$ is naturally
homotopy equivalent to the triple $(G,G_a,G_b)$, and under this
homotopy equivalence the map $f\times_{\f} f_0:T\times_{\f} T_0\to
T\times_{\f} T_0$ induces a representative $f_{\oa}$ of $\oa$ on $G$ that keeps
$G_a$ and $G_b$ invariant and sends $E$ across itself only
once. By induction, there is a homotopy independent of $\oa$ supported in
a small neighborhood of $G_a\cup G_b$
such that $f_{\oa}$ is upper triangular
and such that the restrictions to $G_a$ and $G_b$ satisfy
the conclusions of Theorem \ref{graph kolchin}.

We now claim that if $A$ and $B$ are nonabelian, then the collection
of $f_{\oa}$'s provides the desired lift to $\Q$. We first argue that
if $\oa\in\H$ then $f_{\oa^{-1}}f_{\oa}$ is homotopic to the identity
rel vertices. This map is freely homotopic to the identity and by the
inductive hypothesis it is homotopic rel vertices to a map $g:G\to G$
that is identity on $G_a\cup G_b$ and maps $E$ to a path of the form
$vEu$ where $v$ and $u$ are closed geodesic paths in $G_a$ and $G_b$
respectively. It remains to show that $u$ and $v$ are trivial paths.
Suppose for example that $u$ is nontrivial. We regard the endpoints of $E$
as the basepoints for $G_a$ and $G_b$. Then we may choose closed
paths $\alpha$ and $\beta$ in $G_a$ and $G_b$ so that $\alpha$ does not
commute with $u$ and so that $\beta$ is nontrivial. The closed loop
$E\alpha E^{-1}\beta$ is sent by $g$ to $vEu\alpha u^{-1} E^{-1}v^{-1}\beta$.
Since the two loops are freely homotopic, we conclude that $u$ and $\alpha$
commute, contradicting the choice of $\alpha$. One similarly argues that
$f_{\oa_1}f_{\oa_2}$ is homotopic to $f_{\oa_1\oa_2}$ thus proving the claim
in the case that $A$ and $B$ are nonabelian.

Next suppose that $A$ is abelian and $B$ nonabelian. Then $G_a$ is a
circle with a single edge $\alpha$ and each $f_\oa$ sends $E$ to a
path of the form $\alpha^{m(\oa)}Eu(\oa)$. Define a new filtered graph
$G'=G_b\cup E'$ where $E'$ is a loop based at the basepoint of $G_b$
with the filtration defined so that $E'$ is topmost and the induced
filtration on $G_b$ is unchanged. Define the representative
$f_{\oa}':G'\to G'$ of $\oa$ to agree with $f_\oa$ on $G_b$ and to
send $E'$ to $u(\oa^{-1})E'u(\oa)$. Another way to describe $G'$ is that
it is the result of replacing the ``balloon'' $E\cup \alpha$ with the single
loop $E'$ corresponding to $E^{-1}\alpha E$.
The collection $\{ f_\oa'\}$ forms the desired lift.

If both $A$ and $B$ are abelian, then $\H$ is trivial (by the
preceding argument) and we can take $G$ to be the rose with two petals.

In the case when $T/F_n$ is a circle, i.e. each vertex of $T$ is a translate
of $a$,  we can construct $G$ from $G_a$ by attaching a topmost loop $E$ to
a vertex. The details are entirely analogous to the above discussion of
the nonabelian case and are left to the reader.

From the above discussion it follows that $G$ contains a (unique) maximal
tree such that all edges in the complement are loops, and furthermore
(when $n>1$) each vertex belongs to at least two edges not in the tree.
If $V$ is the number of vertices in $G$, then $G$ has $V-1+n$ edges and
$n\geq 2V$. Thus $V-1+n\leq \frac{3n}{2}-1$ as required.
\qed


\ifx\undefined\bysame
\newcommand{\bysame}{\leavevmode\hbox to3em{\hrulefill}\,}
\fi

\end{document}